\documentclass{lmcs}
\pdfoutput=1

\usepackage{lastpage}
\lmcsdoi{18}{1}{28}
\lmcsheading{}{\pageref{LastPage}}{}{}%
{Apr.~20,~2020}{Feb.~03,~2022}{}

\usepackage[utf8]{inputenc}

\title[Canonicity and homotopy canonicity for cubical type theory]
{Canonicity and homotopy canonicity\texorpdfstring{\\}{} for cubical type theory}

\author[T.~Coquand]{Thierry Coquand\rsuper{a}}
\author[S.~Huber]{Simon Huber\rsuper{a}}
\address{University of Gothenburg, Sweden}
\email{coquand@chalmers.se, simonh@fripost.org}

\author[C.~Sattler]{Christian Sattler\rsuper{b}}
\address{Chalmers University of Technology, Sweden}
\email{sattler@chalmers.se}

\keywords{cubical type theory, univalence, canonicity, sconing, Artin glueing}

\usepackage{mathpartir}
\usepackage{verbatim}
\usepackage[all]{xy}

\usepackage{hyperref}
\usepackage{cleveref}

\usepackage{macros}

\makeatletter
\newcommand{\thmlabel}[1]{%
  \label[\@currenvir]{#1}%
}
\makeatother

\crefname{thm}{Theorem}{Theorems}
\crefname{lem}{Lemma}{Lemmas}
\crefname{prop}{Proposition}{Propositions}
\crefname{rem}{Remark}{Remarks}
\crefname{conj}{Conjecture}{Conjectures}
\crefname{figure}{Figure}{figures}

\crefname{section}{Section}{Sections}
\crefname{subsection}{Subsection}{Subsections}
\crefname{subsubsection}{Subsubsection}{Subsubsections}

\crefname{appendix}{Appendix}{Appendices}

\begin{document}

\begin{abstract}
  Cubical type theory provides a constructive justification of
  homotopy type theory.  A crucial ingredient of cubical type
  theory is a path lifting operation which is explained
  computationally by induction on the type involving several
  non-canonical choices.  We present in this article two canonicity
  results, both proved by a sconing argument:
  a \emph{homotopy canonicity} result, every
  natural number is \emph{path equal} to a numeral, even if we take
  away the equations defining the lifting operation on the type structure, and a
  \emph{canonicity} result, which uses these equations in a crucial way. Both
  proofs are done internally in a presheaf model.
\end{abstract}

\maketitle

\section*{Introduction}

This article is a contribution to the analysis of the computational content of the univalence axiom~\cite{VV} (and higher inductive types).
In previous work~\cite{ABCFHL,BCH,CCHM,CHM,OP}, various {\em presheaf models} of this axiom have been described in a constructive metatheory.
In this formalism, the notion of fibrant type is stated as a refinement of the path lifting operation where one not only provides one of the endpoints but also a partial lift (for a suitable notion of partiality).
This generalized form of path lifting operation is a way to state a homotopy extension property, which was recognized very early (see \eg~\cite{Eilenberg1}) as a key for an abstract development of algebraic topology.
The axiom of univalence is then captured by a suitable \emph{equivalence extension operation} (the ``glueing'' operation), which expresses that we can extend a partially defined equivalence of a given total codomain to a total equivalence.
These presheaf models suggest possible extensions of type theory where we manipulate higher dimensional objects~\cite{ABCFHL,CCHM}.
One can define a notion of reduction and prove canonicity for this extension~\cite{Huber}: any closed term of type $\iN$
(natural number) is convertible to a numeral.
There are however several \emph{non-canonical} choices when defining the path lifting operation
by induction on the type, which produce different notion of convertibility.%
\footnote{For instance, the definition of this operation for ``glue'' types is different in~\cite{CCHM} and~\cite{OP}.}
A natural question is how \emph{essential} these non-canonical choices are: can it be that a closed
term of type $\iN$, defined without use of such non-canonical reduction rules, becomes convertible to $0$ for one choice and $1$ for another?
The main result of this article, the \emph{homotopy canonicity theorem}, implies that this cannot be the case:
the value of a term is independent of these non-canonical choices.
Homotopy canonicity states that, even without providing reduction rules for path lifting operations at type formers,
we \emph{still} have that any closed term of type $\iN$ is \emph{path equal} to a
numeral.
(We cannot hope to have convertibility anymore with these path lifting constants.)
We can then see this numeral as the ``value'' of the given term.

Our proof of homotopy canonicity can be seen as a proof-relevant extension of the \emph{reducibility} or \emph{computability} method, going back to the work of G\"odel~\cite{Dialectica} and Tait~\cite{Tait}.
It is however best expressed in an \emph{algebraic} setting.
We first define a general notion of model, called \emph{cubical category with families}, defined as a category with families~\cite{Dybjer} with certain special operations internal to presheaves over a category $\mathcal{C}$ (such as a cube category) with respect to the parameters of an \emph{interval} $\II$ and a \emph{cofibration classifier} $\FF$.
In this article, we will work with models of the cubical type theory described by~\cite{CCHM,OP}.
However, our methods apply equally well to other versions of cubical type theories that can be presented in a similar setting, for example~\cite{ABCFHL}.

We describe the term model and how to re-interpret the cubical presheaf models as cubical categories with families.
The computability method can then be expressed as a general operation (called ``sconing'') which applied to an arbitrary model $\mathcal{M}$ produces a new model $\MM^*$ with a strict morphism $\MM^* \to \MM$.
Homotopy canonicity is obtained by applying this general operation to the initial model, which we conjecture to be the term model.
This construction associates to a (for simplicity, closed) type $A$ a predicate $A'$ on the closed terms $\verts{A}$ and each closed term $u$ of $A$ a proof $u'$ of $A'\,u$.
The main rules in the closed case are summarized in \cref{fig:sconing}.

We explain next how a similar method can be used to prove {\em canonicity} (or ``strict'' canonicity) when we add computation rules of filling at type formers (using as primitive the operation of {\em composition}).
Here, every closed term of type $\iN$ is (strictly) equal (as opposed to path equal) to a numeral.
This was originally proved by~\cite{Huber}. The main advantage of the present
approach is that we don't need to define an auxiliary reduction relation, and that it is independent
of the exact choice of the equational presentation of cubical type theory.

Some extensions and variations are then described:
\begin{itemize}
\item
Our development extends uniformly to identity types and higher inductive types (using the methods of~\cite{CHM}) (\cref{identity-types,higher-inductive-types}).
\item
Our development applies equally to the case where one treats univalence instead of glue types as primitive (\cref{univalence-as-axiom}).
We expect that a similar sconing argument (glueing along a global sections functor to simplicial sets) works to establish homotopy canonicity for the initial split univalent simplicial tribe in the setting of Joyal~\cite{joyal:tribes}.
\item
Assuming excluded middle, a version of the simplicial set model~\cite{simplicial-set-model} forms an instance of our development, and distributive lattice cubical type theory interprets in it (\cref{simplicial-sets}).
\end{itemize}
Using our technique, one may also reprove canonicity for ordinary Martin-L\"of type theory with inductive families in a reduction-free way.

Shulman~\cite{Shulman} proves homotopy canonicity for homotopy type theory with a truncatedness assumption using the sconing technique.
This proof was one starting point for the present work.

\begin{figure}
\begin{align*}
\iPi(A,B)'(w)
\ &=\ %
\prod{u : \verts{A}} \prod{u' : A'\,u} B'\,u\,u'\,(\iapp(w,u))
\\
\iSigma(A,B)'(w)
\ &=\ %
\sum{u': A'\,(\ifst (w))} B'\,(\ifst(w))\,u'\,(\isnd(w))
\\
\iPath(A,a_0,a_1)'(w)
\ &=\ %
\Path_{\lam{i} A'\,i\,(\ipapp(w,i))}\,a_0'\,a_1'
\\
\iGlue(A,\psi \mapsto (B,w))'(v)
\ &=\ %
\Glue\,\left(A'\,(\iapp(\iunglue,v))\right)\,\left[\psi\mapsto (B'\,v, (w'.1\,v, \ldots))\right]
\end{align*}
\caption{These are the main rules for the computability predicate component in the sconing models for homotopy canonicity and canonicity in the case of the global context. The component relating to fibrancy differs between the two cases.}
\label{fig:sconing}
\end{figure}

\paragraph*{Two models of type theory and sconing}

Since Martin Hofmann's work~\cite{Hofmann97}, it is known how to interpret {\em extensional}
type theory with a hierarchy of universes $\U_0,\U_1,\dots$ in any presheaf model.
As explained in~\cite{CCHM,OP,COQ2018}, in some class of presheaf models, parametrised by
two presheaves $\II$ (representing an abstract interval) and $\FF$ (the cofibration classifier),
it is possible to define, as an internal model inside this presheaf
model, a model of type theory with a hierarchy of universes $\Ufib_0,\Ufib_1,\dots$
satisfying the {\em univalence} axiom.
Both models are carried out in a {\em constructive} metalanguage. In particular, the second
model provides a {\em computational} interpretation of univalence.

This model of univalence is a model of {\em cubical type theory}
where each type has a {\em filling operation}. Univalence is then a {\em theorem} and not
an axiom of cubical type theory. This filling operation is defined by induction
on the type, using a more primitive {\em composition operation}.

The basic scheme for a canonicity proof that we follow here is to associate by induction on a type $A$ a computability predicate $A'$ on the (internal) set $\verts{A}$ of closed elements of this type.
As explained in~\cite{COQ2019}, this so-called sconing interpretation for canonicity goes back to G\"odel's notion of {\em computability} predicates~\cite{Dialectica}, with the crucial feature here that these predicates are now proof-relevant.
This scheme works as well for cubical type theory if we use a metalanguage with an interval object.

We now explain in general terms and by example the differences between the {\em homotopy canonicity} and the {\em canonicity} proofs.
For the {\em homotopy canonicity} proof, we will have $A' : \verts{A} \to \Ufib_n$, while for the {\em canonicity} proof, we will have $A' : \verts{A} \to \U_n$ with a separate component tracking computability of composition.

For the type of natural number $\iN$, for the {\em canonicity} proof, we define $\iN'\,t$ to be $\sum{k:\NN}t =_{\verts{\iN}} \isucc^k(\izero)$, where $t =_{\verts{\iN}} \isucc^k(\izero)$ is (strict) equality, and $\NN$ is the constant presheaf of natural numbers.
In this case, $\iN'$ is not a fibrant family over $\verts{N}$. For {\em homotopy canonicity}, we have to define $\iN'\,t$ as a {\em fibrant} family over $\verts{N}$ (see \cref{sconing:natural-numbers}).%
\footnote{
This fibrant family is not simply obtained by replacing the equality $t =_{\verts{\iN}} \isucc^k(\izero)$ by a path, as this would not model the $\beta$-equality of the eliminator in the successor case.
Instead, we should view $\iN'$ in the case of canonicity as an indexed inductive set and then replace it by a fibrant indexed inductive set for the case of homotopy canonicity.
}

One key step in both arguments is in ensuring that the {\em filling operation}
is a {\em computable} operation.
This is solved in very different ways for the two theories.
For the {\em homotopy canonicity} proof, where $A':\verts{A}\to \Ufib_n$, we can prove
directly that the filling operation is computable without needing information on how the filling operation behaves at individual type formers.
For the {\em canonicity} proof, where $A':\verts{A}\to \U_n$, the filling operation
is defined in terms of a more primitive composition operation and we prove by
induction on the type that this composition operation is computable.

As in~\cite{COQ2019}, we think that this interpretation is best described in an \emph{algebraic} way, using what is essentially a generalized algebraic presentation of type theory.
The difference with~\cite{COQ2019} is that the notion of generalized algebraic theory we are using is now developed internally to a presheaf model with an interval $\II$ and cofibration classifier $\FF$.

\paragraph*{Setting}

We work in a constructive set theory (as presented \eg in~\cite{Aczel}) with a sufficiently long cumulative hierarchy of Grothendieck universes.
However, our constructions are not specific to this setting and can be replayed in other constructive metatheories such as extensional type theory.
In \cref{simplicial-sets}, we assume classical logic for the discussion of models in simplicial sets.

\section{Cubical categories with families}

We first recall the notion of categories with families (cwf)~\cite{Dybjer} equipped with $\Pi$- and $\Sigma$-types, universes, and natural number types.
This notion can be interpreted in any presheaf model.
In that setting, we can consider new operations.
A {\em cubical cwf} will be such a cwf in a presheaf model with extra operations that make use of an interval object~$\II$ and a cofibration classifier~$\FF$ as introduced in~\cite{CHM,OP}.

\subsection{Categories with families}
\label{sec:cwf}

Categories with families form an algebraic notion of model of type theory.
In order to simplify the treatment of universes, we define them in a stratified manner where instead of a single presheaf of types, we specify a filtration of presheaves of ``small'' types.%
\footnote{
We note that this, some might say, non-algebraic aspect of the definition does not interfere with the otherwise algebraic character and that subsets could in principle be replaced by injections.
Indeed, one can even relax the requirement that $\Type_n \to \Type$ is a monomorphism, at the cost of making it more tedious to state coherence of type formers under lifting and level coercion (if desired).
One can also give a version where there is no top-level presheaf of types $\Type$.
None of these variations impact what we do in this article.
}
The length of the filtration is not essential: we have chosen $\omega + 1$ so that we may specify constructions just at the top level.

A \emph{category with families} (cwf) consists of the following data.
\begin{itemize}
\item
We have a category of \emph{contexts} $\Con$ and \emph{substitutions} $\Subst(\Delta, \Gamma)$ from $\Delta$ to $\Gamma$ in $\Con$.
The identity substitution on $\Gamma$ in $\Con$ is written $\id$, and the composition of $\delta$ in $\Subst(\Theta, \Delta)$ and $\sigma$ in $\Subst(\Delta, \Gamma)$ is written $\sigma \delta$.
\item
We have a presheaf $\Type$ of \emph{types} over the category of contexts.
The action of $\sigma$ in $\Subst(\Delta, \Gamma)$ on a type $A$ over $\Gamma$ is written $A \sigma$.
We have a cumulative sequence of subpresheaves $\Type_n$ of \emph{types of level $n$} of $\Type$ where $n$ is a natural number.
\item
We have a presheaf $\Elem$ of \emph{elements} over the category of elements of $\Type$, \ie a set $\Elem(\Gamma, A)$ for $A$ in $\Type(\Gamma)$ with $a \sigma$ in $\Elem(\Delta, A \sigma)$ for $a$ in $\Elem(\Gamma, A)$ and $\sigma$ in $\Subst(\Delta, \Gamma)$ satisfying evident laws.
\item
We have a terminal context $1$, with the unique element of $\Subst(\Gamma, 1)$ written $()$.
\item
Given $A$ in $\Type(\Gamma)$, we have a \emph{context extension} $\Gamma.A$.
There is a \emph{projection} $\pp$ in $\Subst(\Gamma.A, \Gamma)$ and a \emph{generic term} $\qq$ in $\Elem(\Gamma.A, A\pp)$.
Given $\sigma$ in $\Subst(\Delta, \Gamma)$, $A$ in $\Type(\Gamma)$, and $a$ in $\Elem(\Delta, A \sigma)$ we have a \emph{substitution extension} $(\sigma, a)$ in $\Subst(\Delta, \Gamma.A)$.
These operations satisfy $\pp (\sigma, a) = \sigma$, $\qq (\sigma, a) = a$, and $(\pp \sigma, \qq \sigma) = \sigma$.
Thus, every element of $\Subst(\Delta, \Gamma.A)$ is uniquely of the form $(\sigma, a)$ with $\sigma$ and $a$ as above.
\end{itemize}

We introduce some shorthand notation related to substitution.
Given $\sigma$ in $\Subst(\Delta, \Gamma)$ and $A$ in $\Type(\Gamma)$, we write $\sigma^+ = (\sigma \pp, \qq)$ in $\Subst(\Delta.A \sigma, \Gamma.A)$.
Given $a$ in $\Elem(\Gamma, A)$, we write $\subst{a} = (\id, a)$ in $\Subst(\Gamma, \Gamma.A)$.
Thus, given $B$ in $\Type_n(\Gamma.A)$ and $a$ in $\Elem(\Gamma, A)$, we have $B \subst{a}$ in $\Type(\Gamma)$.
Given furthermore $b$ in $\Elem(\Gamma.A, B)$, we have $b \subst{a}$ in $\Elem(\Gamma, B\subst{a})$.
We extend this notation to several arguments: given $a_i$ in $\Elem(\Gamma, A_i)$ for $1 \leq i \leq k$, we write $\subst{a_1, \ldots, a_k}$ for $[a_k] [a_{k-1} \pp] \cdots [a_1 \pp \ldots \pp]$ in
$\Subst(\Gamma, \Gamma.A_1.\ldots.A_k)$.

Note that we could take a different equational presentation. For instance, the presentation in~\cite{ehrhard:phd} takes as primitive the operations $\sigma^+$ and $[u]$ and {\em defines}
then $(\sigma,u)$ as a derived operation $(\sigma,u) = \sigma^+[u]$. It is a strength of the present
approach to canonicity proof to be independent of this choice.

Given a cwf as above, we define what it means to have the following type formers.
In addition to the specified laws, all specified operations are furthermore required to be stable under substitution in the evident manner.
\begin{itemize}
\item
\textbf{Dependent products.}
For $A$ in $\Type(\Gamma)$ and $B$ in $\Type(\Gamma.A)$, we have $\iPi(A, B)$ in $\Type(\Gamma)$, of level $n$ if $A$ and $B$ are.
Given $b$ in $\Elem(\Gamma.A, B)$, we have the \emph{abstraction} $\ilambda(b)$ in $\Elem(\Gamma, \iPi(A, B))$.
Given $c$ in $\Elem(\Gamma, \iPi(A, B))$ and $a$ in $\Elem(\Gamma, A)$, we have the \emph{application} $\iapp(c, a)$ in $\Elem(\Gamma, B \subst{a})$.
These operations satisfy
\begin{align*}
\iapp(\ilambda(b), a) &= b \subst{a}
,\\
\ilambda(\iapp(c \pp, \qq)) &= c
.\end{align*}

%All operations are stable under substitution in that
%\begin{align*}
%\iPi(A, B) \sigma &= \iPi(A \sigma, B \sigma^+)
%,&
%\ilambda(b) \sigma &= \ilambda(b \sigma)
%,&
%\iapp(c, a) \sigma &= \iapp(c \sigma, a \sigma)
%.\end{align*}

Given $A$ and $B$ in $\Type(\Gamma)$ we write $A \to B$ for $\iPi(A, B \pp)$.
\item
\textbf{Dependent sums.}
For $A$ in $\Type(\Gamma)$ and $B$ in $\Type(\Gamma.A)$, we have $\iSigma(A, B)$ in $\Type(\Gamma)$, of level $n$ if $A$ and $B$ are.
Given $a$ in $\Elem(\Gamma, A)$ and $b$ in $\Elem(\Gamma, B \subst{a})$, we have the \emph{pairing} $\ipair(a, b)$ in $\Elem(\Gamma, \iSigma(A, B))$.
Given $c$ in $\Elem(\Gamma, \iSigma(A, B))$, we have the \emph{first projection} $\ifst(c)$ in $\Elem(\Gamma, A)$ and \emph{second projection} $\isnd(c)$ in $\Elem(\Gamma, B \subst{\ifst(c)})$.
These operations satisfy
\begin{align*}
\ifst(\ipair(a, b)) = a
,\\
\isnd(\ipair(a, b)) = b
,\\
\ipair(\ifst(c), \isnd(c)) = c
.\end{align*}
Thus, every element of $\Elem(\Gamma, \iSigma(A, B))$ is uniquely of the form $\ipair(a, b)$ with $a$ and $b$ as above.

%All operations are stable under substitution in that
%\begin{align*}
%\iSigma(A, B) \sigma &= \iSigma(A \sigma, B \sigma^+)
%,&
%\ipair(a, b) \sigma &= \ipair(a \sigma, b \sigma)
%,&
%\ifst(c) \sigma &= \ifst(c \sigma)
%,&
%\isnd(c) \sigma &= \isnd(c \sigma)
%.\end{align*}

Given $A$ and $B$ in $\Type(\Gamma)$ we write $A \times B$ for $\iSigma(A, B \pp)$.
\item
\textbf{Universes.}
We have $\iU_n$ in $\Type_{n+1}(\Gamma)$ and an isomorphism $\Type_n(\Gamma) \cong \Elem(\Gamma, \iU_n)$, naturally in $\Gamma$.%
\footnote{
  This presents Tarski-style universes.
  For Russell-style universes, we would additionally demand that this isomorphism is an identity.
}
%These elements are stable under substitution, \ie $\iU_n \sigma = \iU_n$.
\item
\textbf{Natural numbers.}
We have $\iN$ in $\Type_0(\Gamma)$ with \emph{zero} $\izero$ in $\Elem(\Gamma, \iN)$ and \emph{successor} $\isucc(n)$ in $\Elem(\Gamma, \iN)$ for $n$ in $\Elem(\Gamma, \iN)$.
Given $P$ in $\Type(\Gamma.\iN)$, $z$ in $\Elem(\Gamma, P[\izero])$, $s$ in $\Elem(\Gamma.\iN.P, P (\pp, S(\qq)) \pp)$, and $n : \Elem(\Gamma, \iN)$, we have the \emph{elimination} $\inatrec(P, z, s, n)$ in $\Elem(\Gamma, P[n])$ with
\begin{align*}
\inatrec(P, z, s, \izero) &= z,
\\
\inatrec(P, z, s, \isucc(n)) &= s[n, \inatrec(P, z, s, n)]
.\end{align*}
\end{itemize}
A \emph{structured cwf} is a cwf with type formers as above.

A (strict) morphism $\mathcal{M} \to \mathcal{N}$ of cwfs is defined in the evident manner and consists of a functor $F \co \Con_{\mathcal{M}} \to \Con_{\mathcal{N}}$ and natural transformations $u \co \Type_{\mathcal{M}} \to \Type_{\mathcal{N}} F$ and $v \co \Elem_{\mathcal{M}} \to \Elem_{\mathcal{N}} (F, u)$ such that $v$ restricts to types of level $n$ and the terminal context and context extension is preserved strictly.
A morphism $\mathcal{M} \to \mathcal{N}$ of structured cwfs additionally preserves the operations of the above type formers.
We obtain a category of structured cwfs.

\subsection{Internal language of presheaves}

For the rest of the article, we fix a category $\mathcal{C}$ in the lowest Grothendieck universe.
As in~\cite{ABCFHL,OP,LOPS}, we will use the language of extensional type theory (with subtypes) to describe constructions in the presheaf topos over $\mathcal{C}$.

In the interpretation of this language, a context is a presheaf $A$ over $\mathcal{C}$, a type $B$ over $A$ is a presheaf over the category of elements of $A$, and an element of $B$ is a section.
A \emph{global type} is a type in the global context, \ie a presheaf over $\mathcal{C}$.
Similarly, a \emph{global element} of a global type is a section of that presheaf.

For elements $x$ and $y$ of a type $A$, we have the equality type $x =_A y$, satisfying reflection (we allow ourselves to omit the subscript $A$ if it is evident from the context).
Given a dependent type $B$ over a type $A$, we think of $B$ as a family of types $B\,a$ indexed by elements $a$ of $A$.
We have the usual dependent sum $\sum{a:A} B\,a$ and dependent product $\prod{a:A} B\,a$, with projections of $s : \sum{a:A} B\,a$ written $s.1 : A$ and $s.2 : B\,s.1$, and application of $f : \prod{a:A} B\,a$ to $a : A$ written $f\,a$.
We have also the categorical pairing $\angles{f, g} : X \to \sum{a:A} B\,a$ given $f : X \to A$ and $g : \prod{x:X} B\,(f\,a)$ and other commonly used notations.
The hierarchy of Grothendieck universes in the ambient set theory gives rise to a cumulative hierarchy $\U_0, \U_1, \ldots, \U_\omega$ of universes \`a la Russell.
We model propositions as subtypes of a fixed type $1$ with unique element $\TT$.
This implies that logically equivalent propositions are equal.
We have subuniverses $\Omega_i \subseteq \U_i$ of propositions for $i \in \braces{0, 1 \ldots, \omega}$.

When working in this internal language, we refer to the types as ``sets'' to avoid ambiguity with the types of (internal) cwfs we will be considering.

\subsection{Cubical categories with families}
\label{cubical-cwf}

We now work internally to presheaves over $\mathcal{C}$.
We assume the following:
\begin{itemize}
\item an \emph{interval} $\II : \U_0$ with \emph{endpoints} $0, 1 : \II$,
\item an \emph{cofibration classifier} consisting of $\FF : \U_0$ with a monomorphism $[-] : \FF \to \Omega_0$.%
  \footnote{
    The requirement that $[-]$ is mono is not essential and can be relaxed.
    However, this comes at the cost of making later conditions on $\FF$ more tedious to state.}%
\end{itemize}

As in~\cite{CHM,OP}, a \emph{partial element} of a set $T$ is given by an element $\phi$ in $\FF$ and a function $[\phi] \to T$.
We say that a total element $v$ of $T$ extends such a partial element $\phi, u$ if we have $[\phi] \to u\,\TT = v$.
(Note that the last equation make sense because $[\phi] = \TT$ as soon as $[\phi]$ is inhabited).

Given $A : \II \to \U_\omega$, we write $\hasFill(A)$ for the set of operations taking as inputs $\phi$ in $\FF$, $b \in \braces{0, 1}$, and a partial section $u$ in $\prod{i : \II} [\phi] \vee (i = b) \to A\,i$ and producing an extension of $u$ to a total section in $\prod{i:\II}\, A\,i$.
Given a set $X$ and $Y : X \to \U_\omega$, we write $\Fill(X, Y)$ for the set of \emph{filling structures} on $Y$, producing an element of $\hasFill(Y \circ x)$ for $x$ in $\II \to X$.
Given $s$ in $\Fill(X, Y)$ and $x, \phi, b, u$ as above, we write $s(x, \phi, b, u)$ for the resulting total section in $\prod{i:\II} Y\,(x\,i)$.

We now interpret the definitions of \cref{sec:cwf} in the internal language of the presheaf topos.
A \emph{cubical cwf} is a structured cwf denoted as before that additionally has the following \emph{cubical operations and type formers}.
Again, all specified operations are required to be stable under substitution.
\begin{itemize}
\item
\textbf{Filling operation.}\label{filling}
We have $\ifill$ in $\Fill(\Type(\Gamma), \lam{A} \Elem(\Gamma, A))$ for $\Gamma$ in $\Con$.
Let us spell out stability under substitution: given $A \co \II \to \Type(\Gamma)$, $\phi$ in $\FF$, $b \in \braces{0, 1}$, $u$ in $\prod{i:\II}\, [\phi] \vee (i = b) \to \Elem(\Gamma, A\,i)$, and $\sigma$ in $\Subst(\Delta, \Gamma)$ and $r : \II$, we have
\[
(\ifill(A, \phi, b, u)\,r) \sigma = \ifill(\lam{i} (A\, i) \sigma, \phi, b, \lam{i,x} (u\,i\,x) \sigma)\,r
.\]
\end{itemize}
Note that we do not include computation rules for $\ifill$ at type formers.
This corresponds to our decision to treat $\ifill$ as a non-canonical operation.
\begin{itemize}
\item
\textbf{Dependent path types.}
Given $A$ in $\II \to \Type(\Gamma)$ with $a_b$ in $\Elem(\Gamma, A b)$ for $b \in \braces{0, 1}$, we have $\iPath(A, a_0, a_1)$ in $\Type(\Gamma)$, of level $n$ if $A$ is.
Given $u$ in $\prod{i:\II}\Elem(\Gamma, A i)$, we have the \emph{path abstraction} $\ipabs(u)$ in $\Elem(\Gamma, \iPath(A, u\,0, u\,1))$.
Given $p$ in $\Elem(\Gamma, \iPath(A, a_0, a_1))$ and $i$ in $\II$, we have the \emph{path application} $\ipapp(p, r)$ in $\Elem(\Gamma, A i)$.
These operations satisfy the laws
\begin{align*}
\ipapp(p, b) &= a_b
,\\
\ipapp(\ipabs(u), i) &= u\,i
,\\
\ipabs(\lam{i} \ipapp(p, i)) &= p
.\end{align*}
Thus, every element of $\Elem(\Gamma, \iPath(A, a_0, a_1))$ is uniquely of the form $\ipabs(u)$ with $u$ in $\prod{i:\II}\Elem(\Gamma, A i)$ such that $u\,0 = a_0$ and $u\,1 = a_1$.
%All operations are stable under substitution in that
%\begin{align*}
%\iPath(A, a_0, a_1) \sigma &= \iPath(A \sigma, a_0 \sigma, a_1 \sigma)
%,&
%\ipabs(u) \sigma &= \ipabs(\lam{i} (u i) \sigma)
%,&
%\ipapp(p, i) \sigma &= \ipapp(p \sigma, r).
%\end{align*}
\end{itemize}
Using path types, we define $\iisContr(A)$ in $\Type(\Gamma)$ for $A$ in $\Type(\Gamma)$ as well as $\iisEquiv$ in $\Type(\Gamma.A \to B)$ and $\iEquiv$ in $\Type(\Gamma)$ for $A, B$ in $\Type(\Gamma)$ as in~\cite{CCHM}.
(We use a subscript here and for some other notions to distinguish them from analogous notions defined later in a different setting in \cref{developments}.)
These notions are used in the following type former, which extends any partially defined equivalence (given total codomain) to a totally defined function.
\begin{itemize}
\item
\textbf{Glue types.}
Given $A$ in $\Type(\Gamma)$, $\phi$ in $\FF$, $T$ in $[\phi] \to \Type(\Gamma)$, and
\[
e : [\phi] \to \Elem(\Gamma, \iEquiv(T\,\TT, A))
,\]
we have the \emph{glueing} $\iGlue(A, \phi, T, e)$ in $\Type(\Gamma)$, equal to $T$ on $[\phi]$ and of level $n$ if $A$ and $T$ are.
We have $\iunglue$ in $\Elem(\Gamma, \iGlue(A, \phi, T, e) \to A)$ such that $\iunglue = \ifst(e)\,\TT$ on $[\phi]$.
Given $a$ in $\Elem(\Gamma, A)$ and $t$ in $[\phi] \to \Elem(\Gamma, T)$ such that $\iapp(\ifst(e)\,\TT, t\,\TT) = a$ on $[\phi]$, we have $\iglue(a, t)$ in $\Elem(\Gamma, \iGlue(A, \phi, T, e))$ equal to $t$ on $[\phi]$.
These operations satisfy
\begin{align*}
\iapp(\iunglue, \iglue(a, t)) &= a
,\\
\iglue(\iapp(\iunglue, u), \lam{x} u) &= u
.\end{align*}
Thus, every element of $\Elem(\Gamma, \iGlue(A, \phi, T, e))$ is uniquely of the form $\iglue(a, t)$ with $a$ and $t$ as above.
%All operations are stable under substitution in that
%\begin{align*}
%\iGlue(A, \phi, \angles{T, e}) \sigma = \iGlue(A \sigma, \phi, \lam{x} ((T\,x) \sigma, (e\,x) \sigma))
%,\\
%\begin{aligned}
%\iunglue \sigma = \iunglue
%&&,
%\iglue(a, t) = \iglue(a \sigma, t \sigma)
%\end{aligned}
%.\end{align*}
\end{itemize}

The notion of morphism of structured cwfs lifts to an evident notion of morphism of cubical cwfs.
We obtain, internally to presheaves over $\mathcal{C}$, a category of cubical cwfs.
We now lift this category of cubical cwfs from the internal language to the ambient theory by interpreting it in the global context: externally, a cubical cwf (relative to the chosen base category $\mathcal{C}$, interval $\II$, and cofibration classifier $\FF$) consists of a presheaf $\Con$ over $\mathcal{C}$, a presheaf $\Type$ over the category of elements of $\Con$, \etc.

\begin{rem} \thmlabel{univalent-type-theory-prop-J}
Fix a cubical cwf as above.
Assume that $\II$ has a connection algebra structure and that $\FF$ forms a sublattice of $\Omega_0$ that contains the interval endpoint inclusions.
As in~\cite{CCHM}, it is then possible in the above context of the glue type former to construct an element of $\Elem(\Gamma, \iisEquiv[\iunglue])$.
%Thus, glue types show how to extend any partially defined equivalence (given total codomain) to a totally defined equivalence.
From this, one derives an element of $\Elem(\Gamma, \iUnivalence_n)$ where
\[
\iUnivalence_n = \iPi(\iU_n, \iisContr(\iSigma(\iU_n, \iEquiv(\qq, \qq \pp))))
\]
for $n \geq 0$, \ie univalence is provable.
One may also show that the path type applied to constant families $\II \to \Type(\Gamma)$ interprets the rules of identity types of Martin-L\"{o}f with the computation rule for the eliminator $\iJ$ replaced by a propositional equality.
Thus, we obtain an interpretation of univalent type theory with identity types with propositional computation in any cubical cwf.
\end{rem}

\subsection{Computational cubical categories with families}
\label{computational-cubical-cwf}

In this subsection, we consider a variation of the notion of cubical categories with
families where we replace the filling operation by a composition operation,
and where we add {\em computation} rules for this composition operation. This version
is the one used for (strict) canonicity in \cref{strict-canonicity}.
The computation rules are needed since
the proof of canonicity follows closely the constructive justification of cubical type theory.
For this justification, we also have to replace the {\em filling}
operation by a {\em composition operation}. We show then that we can define
a filling operation from a composition operation and we define the composition
operation on types structurally. Since the canonicity argument, like
the one in~\cite{COQ2019}, follows closely the structure of the constructive
justification of the model of univalence, we need to start from the composition
operation instead.

In order to simplify the notations, we assume here that the interval $\II$
also has a reversal operation, like in~\cite{CCHM,CHM}.
This assumption is not necessary (for instance, as noted in~\cite{CCHM}, and indeed as we did in \cref{cubical-cwf},
we can avoid the reverse operation at the cost of carrying around an external boolean
parameter) but it simplifies the presentation slightly.

Given $A : \II \to \U_\omega$, we write $\hasComp(A)$ for the set of operations taking as
inputs $\phi$ in $\FF$ with a partial section $u$ in $\prod{i : \II} [\phi] \vee (i = 0) \to A\,i$
and producing an element in $A\,1$ which is equal to $u\,1\,\TT$ on $\phi$.
Given a set $X$ and $Y : X \to \U_\omega$, we write $\Comp(X, Y)$ for
the set of \emph{composition  structures} on $Y$, producing an element of $\hasComp(Y \circ x)$
for $x$ in $\II \to X$.
Given $s$ in $\Comp(X, Y)$ and $x, \phi, u$ as above, we write $s(x, \phi, u)$ for
the resulting element in $Y\,(x\,1)$.
%\note{What is $b$? Introduce $b$-sided composition and eliminate reverse?
%Note that the versions with and without a parameter $b$ are used inconsistently in the following.}

%% We consider the constant presheaf $\NN$ of natural numbers.
%% One important property is that any function $\II\to N$ is {\em constant}.
%% This follows from the fact that $\II$ is tiny \note{We have not assumed this here. This assumption is introduced in \cref{developments}}.
%% As a consequence, $\NN$ (seen as a family over the terminal set $1$) has a unique composition structure (given by $s(\psi,u) = u_0\,\TT$ for any family $u_i : N$ over $\II$ defined on $(i = 0)\vee\psi$.

We now change the definition of cubical cwf in two ways to obtain our notion of \emph{computational cubical cwf}.

First, we replace the {\em filling} operation by a {\em composition operation}.
We have $\icomp$ in $\Comp(\Type(\Gamma), \lam{A} \Elem(\Gamma, A))$ for $\Gamma$ in $\Con$
together with stability under substitution: given $A \co \II \to \Type(\Gamma)$, $\phi$ in $\FF$,
$u$ in $\prod{i:\II}\, [\phi] \vee (i = 0) \to \Elem(\Gamma, A\,i)$, and $\sigma$
in $\Subst(\Delta, \Gamma)$ and $r : \II$, we have
\[
\icomp(A, \phi, u) \sigma = \icomp(\lam{i} (A\, i) \sigma, \phi, \lam{i,x} (u\,i\,x) \sigma)
.\]

The filling operation is now a {\em derived} operation. We define
\[
\ifill(A,\phi,u)\,r = \icomp(A_r,\phi,u_r)
\]
where $A_r\,i = A\,(i\wedge r)$ and $u_r\,i = u\,(i\wedge r)$.

Second, we add suitable {\em computation rules} (equalities) for this composition operation, structurally over types, following
the computation rules in~\cite{CCHM,CHM}. We give the details here for two representative examples.
\begin{itemize}
\item
For dependent sums, we add the computation rule
\[
\icomp(\lam{i} \iSigma(A\,i,B\,i),\psi,w) = \ipair(u,v)
\]
where $u = \tilde{u}\,1$ and $v = \icomp(\lam{i} (B\,i)[\tilde{u}\,i],\psi,\lam{i,x} \isnd(w\,i\,x))$ using
\[
\tilde{u} = \ifill(A,\psi,\lam{i,x} \ifst(w\,i\,x))
.\]
\item
For natural numbers, we add the computation rules
\begin{align*}
\icomp(\lam{i} \iN,\psi,\lam{i,x} \izero) &= \izero
\\
\icomp(\lam{i} \iN,\psi,\lam{i,x} \isucc (v\,i\,x)) &= \isucc (\icomp(\lam{i} \iN,\psi,b,v))
\rlap{\text{.}}\end{align*}
\end{itemize}

\section{Two examples of cubical cwfs}

In this section we give two examples of cubical cwfs: a term model and
a particular cubical cwfs formulated in a constructive metatheory, the
latter with extra assumptions on $\II$ and $\FF$.

\subsection{Term model}

We sketch how to give a cubical cwf $\Syn$ built from syntax, and
refer the reader to \cref{term-model-rules} for more details.  All our
judgments will be indexed by an object $X$ of $\mathcal{C}$ and given
a judgment $\Gamma \der_X \JJ$ and $f \co Y \to X$ in $\mathcal{C}$ we
get $\Gamma f \der_Y \JJ f$.  Here, $f$ acts on expressions as an
implicit substitution, while for substitutions on object variables we
will use explicit substitutions.

The forms of judgment are:
\begin{mathpar}
   \Gamma \der_X \and %
   \Gamma \der_X A \and %
   \Gamma \der_X A = B \and %
   \Gamma \der_X t : A \and %
   \Gamma \der_X t = u : A \and %
   \sigma \co \Delta \to_X \Gamma
\end{mathpar}
The main rules are given in the appendix.  This then induces a cubical
cwf $\Syn$ by taking, say, the presheaf of contexts at stage $X$ to be
equivalence classes of $\Gamma$ for $\Gamma \der_X$ where the
equivalence relation is judgmental equality.

Some rules are a priori infinitary, but in some cases (such as the one
considered in~\cite{CCHM}) it is possible to present the rules in a
finitary way.

\medskip

This formal system expresses the laws of cubical cwfs in rule form. It defines
the {\em term model}.
Following~\cite{Streicher91,palmgren-vickers:partial-horn-logic} developed in an
intuitionistic framework, we conjecture that this can be interpreted in an
arbitrary cubical cwf in the usual way:

\begin{conj} \thmlabel{thm:initial-ccwf}
With chosen parameters $\mathcal{C}, \II, \FF$, the cubical cwf $\Syn$ is initial in the category of cubical cwfs.
\end{conj}

However, our canonicity result is orthogonal to this conjecture: It is a result about the initial model, without need for an explicit description of this model as a term model.

\subsection{Developments in presheaves over \texorpdfstring{$\mathcal{C}$}{C}}
\label{developments}

We now assume that $\II$ and $\FF$ satisfy the axioms presented in~\cite{OP,COQ2018}.
We briefly recall them for the reader's convenience.
The subobject $\FF$ of $\Omega_0$ should define a {\em dominance} and be closed under disjunction.
The subobject classified by the map $[-] \colon \FF \to \Omega_0$ should be levelwise decidable.
The interval $\II$ should have two distinct global elements $0$ and $1$ and connections.
The interval endpoint inclusions should be cofibrations (\ie, the equalities to $0$ and $1$ are coded by elements of $\FF$) and cofibrations should be closed under universal quantification over $\II$.
Finally, the interval $\II$ should be {\em tiny}, \ie, the exponential functor $(-)^{\II}$ should have a right adjoint $R$.%
\footnote{This is not part of the axioms in~\cite{OP}, but it implies connectivity of $\II$, the first axiom in~\cite{OP}, since left adjoints preserve colimits.}
This is for example the case if $\mathcal{C}$ has finite products and $\II$ is representable.

Most of the reasoning will be done in the internal language of the presheaf topos.
At certain points however, we need to consider the set of global sections of a global type $F$; we denote this by $\Box F$.
We stress that statements involving $\Box$ are external, not to be interpreted in the internal language.
Crucially, the adjunction $(-)^\II \dashv R$ cannot be made internal~\cite{LOPS}.

\newcommand{\C}{\mathcal{C}}
\newcommand{\PC}{\widehat{\C}}

We write $\PC$ for the category of presheaves over $\C$.
The right adjoint $R$ is determined by an isomorphism
\[
\PC(A, R X) \simeq \PC(A^\II, X)
\]
natural in $A$ and $X$.
Using cocontinuity in $X$, we may equivalently restrict to $A = \yon(I)$ where $\yon$ denotes the Yoneda embedding and $I$ is in $\C$.
Then the isomorphism becomes
\[
(R X)(I) \simeq \PC(\yon(I)^\II, X)
\]
natural in $I$ and $X$.
We may modify the given right adjoint $R$ so that this isomorphism becomes an equality.
By our smallness assumptions on $\C$ and $\II$, we have that $y(I)^\II$ lives in the lowest Grothendieck universe in our hierarchy.
It follows that $R$ restricts to an operation on $\U_n$ for $n \geq 0$.%
\footnote{Without our modification of $R$, this would only be true up to isomorphism.} %

Pseudofunctorially in a presheaf $A$, the adjunction $(-)^\II \dashv R$ descends to an adjunction between categories of families over $A$ and $A^\II$.
We record what we need from this in the rest of our development.

\begin{lem} \thmlabel{adjunction-slice}
Let $A$ be a global set and $B$ a global family over $A^\II$.
Then we have a global family $B_\II$ over $A$ with a bijection of global elements
\[
\Box (\prod{A'} B_\II \circ f) \simeq \Box (\prod{(A')^\II} B \circ f^\II)
\]
natural in global $f : A' \to A$.

The construction $(-)_\II$ may be chosen so that:
\begin{enumerate}
\item \label{adjunction-slice:size}
if $B$ is valued in $\U_n$ for $n \geq 0$, then so is $B_\II$,
\item \label{adjunction-slice:substitution}
the induced isomorphism $(B \circ f^\II)_\II \simeq B_\II \circ f$ is an identity.
\end{enumerate}
\end{lem}

\begin{proof}
Let $\eta_A \co A \to R(A^\II)$ be the unit of the adjunction at $A$.
We define $B_\II(a)$ as the fiber of $R((-).1) : R(\sum{A^\II} B) \to R(A^\II)$ over $\eta_A(a)$.
Global sections of $\prod{A'} B_\II \circ f^\II$ are dotted maps making the following diagram commute:
\[
\xymatrix{
&&
  R(\sum{A^\II} B)
  \ar[d]^{R((-).1)}
\\
  A'
  \ar[r]_{f}
  \ar@{.>}[urr]
&
  A
  \ar[r]_{\eta_A}
&
  R(A^\II)
\rlap{.}}
\]
Global sections of $\prod{(A')^\II} B \circ f$ are dotted maps making the following diagram commute:
\[
\xymatrix{
&&
  \sum{A^\II} B
  \ar[d]^{(-).1}
\\
  (A')^\II
  \ar[rr]_{f^\II}
  \ar@{.>}[urr]
&&
  A^\II
\rlap{.}}
\]
Under transposition of the adjunction, the two are in bijection, naturally in $A'$.

Recall the equivalence between maps into a presheaf and families over that presheaf.
Under this equivalence, we can regard $(-)^\II$ as a functor from global families over $A$ to global families over $A^\II$.
The above discussion then shows that $(-)_\II$ is a right adjoint of $(-)^\II$.

Let examine the values of the presheaf $B_\II$ over the category of elements of $A$.
By Yoneda and the natural bijection we have just verified, $B_\II(I, a)$ is naturally isomorphic to the set of sections of the restriction of $B$ along $a^\II \colon \yon(I)^\II \to A^\II$.
By our smallness assumptions on $\C$ and $\II$, this is in $\U_n$ if $B$ is valued in $\U_n$ (irrespective of the size of $A$).
As in our discussion on size preservation of $R$, we may modify the definition of $(-)_\II$ so that the above isomorphism becomes an identity.
This validates~\eqref{adjunction-slice:size} and~\eqref{adjunction-slice:substitution}.
\end{proof}

In the above statement, the given bijection may be reduced to the case where $f$ is an identity: $\Box (\prod{A^\II} B) \simeq \Box (\prod{A} B_{\II})$.
The cost to pay is that the isomorphism $B_\II \circ f \simeq (B \circ f^\II)_\II$ with appropriate coherence becomes primitive (non-derived) data.

We are going to apply \cref{adjunction-slice} in two different instances.
The first instance occurs in the following subsection and is used to build internal universes of fibrant sets, which will be needed to prove homotopy canonicity.
The second (and more complex) instance occurs in in \cref{canonicity-sconing-cwf} and is used for interpreting types in the sconing model in the proof of canonicity.

\subsubsection{Fibrant presheaves}
\label{fibrant}

Recall the global family $\hasFill : \U_\omega^\II \to \U_\omega$ from \cref{cubical-cwf}.
Applying \cref{adjunction-slice} to $\hasFill$, we obtain global $\CC : \U_\omega \to \U_\omega$
such that naturally in a global set $X$ with global $Y : X \to \U_\omega$, global elements of $\prod{x : X^\II} \hasFill(Y \circ x)$ are in bijection with global elements of $\prod{x : X} \CC(Y\,x)$.
Given a global set $X$ and global $Y : X \to \U_\omega$, we thus have a logical equivalence (maps back and forth)
\begin{equation} \label{C-vs-fill-external}
\Box \Fill(X, Y) \longleftrightarrow \Box \prod{x:X}\CC(Y\,x)
\end{equation}
natural in $X$.%
\footnote{We record only the logical equivalence instead of an isomorphism so that it will be easier to apply our constructions in situations where the right adjoint $R$ fails to exist such as \cref{simplicial-sets}.
  Naturality is only used at one point below, for the forward map, to construct suitable elements of $C$ applied to glueings.}

Note that $\CC$ descends to $\CC : \U_n \to \U_n$ for $n \geq 0$.
We write $\Ufib_i = \sum{A:\U_i} \CC(A)$ for $i \in \braces{0, 1, \ldots, \omega}$; we call $\Ufib_i$ a universe of \emph{fibrant sets}.
Now set $X = \Ufib_\omega$ and $Y(A, c) = A$ in~\eqref{C-vs-fill-external}.
We trivially have $\Box \prod{x:X}\CC(Y\,x)$, thus get
\begin{equation} \label{U-fib-has-fill}
\fillop : \Fill(\Ufib_\omega, \lam{(A, c)} A)
.\end{equation}
This is essentially the counit of the adjunction defining $\CC$.
Note that~\cite{LOPS} use modal extensions of type theory to perform this reasoning internal to presheaves over $\mathcal{C}$.

\begin{rem} \thmlabel{C-to-fill}
Internally, a map $\Fill(X, Y) \to \prod{x:X}\CC(Y\,x)$ does not generally exist for a set $X$ and $Y : X \to \U_\omega$ as for $X = 1$ one would derive a filling structure for any ``homogeneously fibrant'' set, which is impossible (see~\cite[Remark~5.9]{OP}).
%\note{What about a counterexample where $X$ and $Y$ are global?}
However, from~\eqref{U-fib-has-fill} we get a map $\prod{x:X} \CC(Y\,x) \to \Fill(X, Y)$ natural in $X$ using closure of filling structures under substitution (see below).
\end{rem}

More examples of the interplay between \emph{internal} and \emph{external} reasoning involving elements of $\CC(X)$ will occur in \cref{general-constructions} when we reason that closure of $\Fill$ under various type formers, proven internally, transfers to corresponding closure properties of $\CC$, proven externally.

\subsubsection{Some general constructions}
\label{general-constructions}

We recall some constructions of~\cite{CCHM,OP} in the internal language.
\begin{itemize}
\item
Given $A : \II \to \U_\omega$ and $a_b : A_b$ for $b \in \braces{0, 1}$, \emph{dependent paths} $\Path_A\,a_0\,a_1$ are the set of maps $p : \prod{i : \II} A\,i$ such that $p\,0 = a_0$ and $p\,1 = a_1$.
We use the same notation for non-dependent paths.
\item
For $A : \U_\omega$, we have a set $\isContr(A)$ of witnesses of \emph{contractibility}, defined using paths.
\item
Given $A, B : \U_\omega$ with $f : A \to B$, we have the set $\isEquiv(f)$ with elements witnessing that $f$ is an \emph{equivalence}, defined using contractibility of homotopy fibers.
We write $\Equiv(A, B) = \sum{f : A \to B} \isEquiv(f)$.
\item
Given $A : \U_\omega$, $\phi : \FF$, $B : [\phi] \to \U_\omega$, and $e : [\phi] \to (B\,\TT \to A)$, the \emph{glueing} $\Glue\,A\,[\phi \mapsto (B, e)]$ consists of elements $\glue\,a\,[\phi \mapsto b]$ with $a : A$ and $b : [\phi] \to B$ such that $e.1\,(b\,\TT) = a$ on $[\phi]$ and is defined in such a way that
\begin{align*}
\Glue\,A\,[\phi \mapsto (B, e)] &= T\,\TT,
\\
\glue\,a\,[\phi \mapsto b] &= b\,\TT
\end{align*}
on $[\phi]$.
We have a projection $\unglue : \Glue\,A\,[\phi \mapsto (B, e)] \to A$.
\end{itemize}
These operations are valued in $\U_n$ if their inputs are.
We further recall from~\cite{CCHM,OP} basic facts about filling structures in the internal language.
\begin{itemize}
\item
Filling structures are closed under substitution: given $f : X' \to X$ and $Y \co X \to \U_\omega$, any element of $\Fill(X, Y)$ induces an element of $\Fill(X', Y \circ f)$, naturally in $X'$.
\item
Filling structures are closed under exponentiation: given sets $S, X$ and $Y \co X \to \U_\omega$, any element of $\Fill(X, Y)$ induces an element of
\[
\Fill(X^S, \lam{x} \prod{s:S} Y(x\,s))
,\]
naturally in $S$.
\item
Filling structures are closed under $\Pi, \Sigma, \Path$.
For example, for dependent products, this means the following.
Given $A : \Gamma \to \U_\omega$ with $\Fill(\Gamma, A)$ and $B : \prod{\rho : \Gamma} A\,\rho \to \U_\omega$ with $\Fill(\sum{\rho : \Gamma} A\,\rho, \lam{(\rho, a)} B\,\rho\,a)$, we have
\[
\Fill(\Gamma, \lam{\rho : \Gamma} \prod{a : A\,\rho} B\,\rho\,a)
.\]
\item
The $\Glue$ set former preserves filling structures with equivalences.
By this, we mean the following.
Let $A : \Gamma \to \U_\omega$ and $\phi : \FF$.
For $\rho$ in $\Gamma$ and $x$ in $[\phi]$, let $B\,\rho\,x : \U_\omega$ with a map $e\,\rho\,x : A\,\rho \to B\,\rho\,x$.
Assume $\Fill(\Gamma, A)$ and $\Fill(\sum{\rho : \Gamma} [\phi], \lam{(\rho, x)} B\,\rho\,\TT)$ and that $e\,\rho\,\TT$ is an equivalence for $\rho : \Gamma$ on $[\phi]$.
Then we have
\[
\Fill(\Gamma, \lam{\rho} \Glue\,(A\,\rho))\,[\phi\,\rho \mapsto (B\,\rho\,\TT, e\,\rho\,\TT)]
\]
and the map $\unglue\,\rho$ from the glue set to $A\,\rho$ is an equivalence for $\rho$ as above.
\end{itemize}
All of the above closure observations satisfy naturality under substitution.

Above, we have recorded closure of $\Fill$ under various set formers.
From this, we use external reasoning to deduce the corresponding closure properties for $\CC$.
Specifically, we have that $\CC$ is closed under $\Pi, \Sigma, \Path, \Glue$ (adding equivalence data in the case of $\Glue$), and that $\CC(A)$ implies $\CC(A^S)$ for $A, S : \U_\omega$.%
\footnote{Note that naturality in $S$ of the latter operation is used in substitutional stability of universes in the sconing in \cref{sconing}.} %

We explain how this works in the example case of $\Pi$.

Given $(A, c_A)$ in $\Ufib_\omega$ and $\angles{B, c_B} : A \to \Ufib_\omega$, we wish to show $C(\prod{A} B)$.
We set
\[
\Delta = \sum{(A, c_A) : \Ufib_\omega} A \to \Ufib_\omega
\]
for the ``generic context'' of the closure statement.
Then the goal is a global element of
\[
\prod{((A, c_A), \angles{B, c_B}) : \Delta} C(\prod{A} B)
.\]
By~\eqref{C-vs-fill-external}, this amounts to a global element of
\[
\Fill(\Delta, \lam{((A, c_A), \angles{B, c_B})} \prod{A} B)
.\]
Now we reason internally.
Since $\Fill$ is closed under dependent products, the goal reduces to
\begin{align*}
&\Fill(\Delta, \lam{((A, c_A), \angles{B, c_B})} A)
,\\
&\Fill(\sum{((A, c_A), \angles{B, c_B}) : \Delta} A, \lam{(((A, c_A), \angles{B, c_B}), a)} B\,a)
.\end{align*}
Elements of these are given by \cref{C-to-fill} since all families here are valued in fibrant sets (as witnessed by the components $c_A$ and $c_B$).

Note that in the case of $\Glue$ with $(A, c) : \Ufib_\omega$, $\phi : \FF$, $\angles{B, d} : [\phi] \to \Ufib_\omega$, and $e : [\phi] \to \Equiv(B\,\TT, A)$, naturality of the forward map of~\eqref{C-vs-fill-external} is needed to see that the element $c : \CC(\Glue\,A\,[\phi \mapsto (B, e)])$ constructed in the same fashion as above for dependent products equals $d\,\TT : C(B\,\TT)$ on $[\phi]$.
%\note{Should this be explained further?}
% Justification for the above

\begin{comment}
  We apply~\eqref{C-vs-fill-external} with a suitable ``generic context''
  $X = \sum{(A, c) : \Ufib_\omega}\sum{\phi : \F}\sum{T : [\phi] \to \Ufib_\omega}
  [\phi] \to \Equiv(T\,\TT, A)$ and $Y(A, \phi, T, e) = \Glue\A\,[\phi \mapsto (T, e)]$.

  X = \sum{A : \Ufib_\omega}\sum{\phi : \FF}\sum{B : [\phi] \to \Ufib_\omega}
  \sum{e : [\phi] \to \Equiv(B\,\TT, A)}
  Y(A, \phi, B, e) = \Glue\,A\,(\phi \mapsto (B, e))

X' = \sum{A : \Ufib_\omega}\sum{B : \Ufib_\omega}\sum{e : \Equiv(B\,\TT, A)}
Y'(A, B, e) = B

X'' = \Ufib_\omega
Y''(B) = B

f : X' \to X
f(A, B, e) = (A, \top, \lam{x} B, \lam{x} e)

g : X' \to X''
g(A, B, e) = B

c : Fill(X, Y)
c'' : Fill(X'', Y'')

the restrictions of c along f and c'' along g to Fill(X', Y') agree
\end{comment}

As in~\cite{CCHM,OP,LOPS}, glueing shows $\Fill(1, \Ufib_n)$ for $n \geq 0$.
Using~\eqref{C-vs-fill-external}, we conclude $\CC(\Ufib_n)$.

Let $\mathbb{N}$ denote the natural number object in presheaves over $\mathcal{C}$, the constant presheaf with value the natural numbers.
From~\cite{CCHM,OP}, we have $\Fill(1, \mathbb{N})$.
Using~\eqref{C-vs-fill-external}, we conclude $\CC(\mathbb{N})$.

We justify fibrant indexed inductive sets in \cref{indexed-inductive-sets}.

\subsection{Standard model}
\label{standard-model}

Making the same assumptions on $\mathcal{C}, \II, \FF$ as in \cref{developments}, we can now specify the standard model $\Std$ of cubical type theory in the sense of the current article as a cubical cwf (with respect to parameters $\mathcal{C}, \II, \FF$) purely using the internal language of the presheaf topos.
The cwf is induced by the family over $\Ufib_\omega$ given by the first projection as follows.
\begin{itemize}
\item
The category of contexts is $\U_{\omega}$, with $\Subst(\Delta, \Gamma)$ the functions from $\Delta$ to $\Gamma$.
\item
The types over $\Gamma$ are maps from $\Gamma$ to $\Ufib_\omega$; a type $\angles{A, p}$ is of level $n$ if $A$ is in $\Gamma \to \U_n$.
This is clearly functorial in $\Gamma$.
\item
The elements of $\angles{A, p} : \Gamma \to \Ufib_\omega$ are $\prod{\rho:\Gamma}A\,\rho$.
This is clearly functorial in $\Gamma$.
\item
The terminal context is given by $1$.
\item
  The context extension of $\Gamma$ by $\angles{A, p}$ is given by
  $\sum{\rho:\Gamma} A \, \rho$, with $\pp, \qq$ given by projections and substitution extension given by pairing.
\end{itemize}
We briefly go through the necessary type formers and operations, omitting evident details.
Whenever we mention an induced witness of fibrancy, this refers to the observations recorded in \cref{general-constructions}.
\begin{itemize}
\item
The dependent product of $\angles{A, c} : \Gamma \to \Ufib_\omega$ and $\angles{B, d} : \sum{\rho : \Gamma} A\,\rho \to \Ufib_\omega$ is
\[
\angles{\lam{\rho} \prod{a : A\,\rho} B(\rho, a), e}
\]
where $e\,\rho : \CC(\prod{a : A\,\rho} B(\rho, a))$ is induced by $c\,\rho : \CC(A\,\rho)$ and $d\,\rho\,a : \CC(B(\rho, a))$ for $a : A$.
\item
The dependent sum of $\angles{A, c} : \Gamma \to \Ufib_\omega$ and $\angles{B, d} : \sum{\rho : \Gamma} A\,\rho \to \Ufib_\omega$ is
\[
\angles{\lam{\rho} \sum{a : A\,\rho} B(\rho, a), e}
\]
where $e$ is induced by $c$ and $d$.
\item
The universe $\iU_n : \Gamma \to \Ufib_{n+1}$ is constantly
\[
(\Ufib_n, c)
\]
with $c : \CC(\Ufib_n)$ as recorded before.
According to our definition of the types in $\Std$, this universe is actually Russell-style, \ie, the evident isomorphism $\Type_n(\Gamma) \cong \Elem(\Gamma, \iU_n)$ is an identity.
\item
The natural number type $\iN : \Gamma \to \Ufib_0$ is constantly
\[
(\mathbb{N}, c)
\]
with $c : \CC(\mathbb{N})$ as recorded before.
The zero and successor constructors and the eliminator are given by the corresponding features of the natural number object $\mathbb{N}$.
\end{itemize}
We now turn to the cubical aspects.
\begin{itemize}
\item
The filling operation
\[
\ifill : \Fill(\Gamma \to \Ufib_\omega, \lam{\angles{A, p}} \prod{\rho : \Gamma} A\,\rho)
\]
is derived from~\eqref{U-fib-has-fill} by closure of filling structures under exponentiation.
\item
Given $\angles{A, c} : \II \to \Gamma \to \sum{A : \U_\omega} \CC(A)$ and $a_b : \prod{\rho : \Gamma} A\,b\,\rho$ for $b \in \braces{0, 1}$, we define $\iPath(A, a_0, a_1) : \Gamma \to \sum{A : \U_\omega} \CC(A)$ as
\[
\angles{\prod{\rho : \Gamma} \Path_{\lam{i} A\,i\,\rho}\,c_0\,c_1), d}
\]
where $d\,\rho : \CC(\Path_{\lam{i} A\,i\,\rho}\,c_0\,c_1))$ is induced by $\lam{i} c\,i\,\rho$.
Path abstraction and application operations are defined from those of $\Path$.
\end{itemize}
Before defining glue types, we note that the notions $\iisContr$ and $\iisEquiv$ in the cubical cwf we are defining correspond to the notions $\isContr$ and $\isEquiv$.
For example, given a type $A : \Gamma \to \Ufib_\omega$, then the elements of $\iisContr(A)$, given by
$\prod{\rho : \Gamma}\,\iisContr(A).1\,\rho$, are in bijection with $\prod{\rho : \Gamma} \isContr(A.1\,\rho)$ naturally in $\Gamma$.
\begin{itemize}
\item
Given $\angles{A, c} : \Gamma \to \Ufib_\omega$, $\phi : \FF$, $\angles{T, d} : [\phi] \to \Gamma \to \Ufib_\omega$, and $e : [\phi] \to \iEquiv(T\,\TT, A)$, we define $\iGlue(\angles{A, c}, \phi, \angles{T, d}, e) : \Gamma \to \Ufib_\omega$ as
\[
\lam{\rho} (\Glue\,(A\,\rho)\,[\phi \mapsto (T\,\TT\,\rho, (e'\,\TT\,\rho).1)], q\,\rho)
\]
where $e'\,\TT\,\rho : \Equiv(T\,\TT\,\rho, A\,\rho)$ is induced by $e\,\TT\,\rho$ and $q\,\rho$ is induced by $c\,\rho$ and $\lam{x} d\,x\,\rho$ and $\lam{x} (e'\,\TT\,\rho).2$.
\end{itemize}

We have thus verified the following statement.

\begin{thm}
Assuming the parameters $\mathcal{C}, \II, \FF$ satisfy the assumptions of \cref{developments}, the standard model $\Std$ forms a cubical cwf.
\end{thm}

\section{Sconing}
\label{sconing}

We make the same assumptions on our parameters $\mathcal{C}, \II, \FF$ as in \cref{developments}.
In the global context, let $\MM$ be a cubical cwf (with respect to these parameters) denoted $\Con, \Subst, \ldots$ as in \cref{cubical-cwf}.
We assume that $\MM$ is \emph{size-compatible} with the standard model, by which we mean $\Subst(\Delta, \Gamma) : \U_\omega$ for all $\Gamma, \Delta$ and $\Elem(\Gamma, A) : \U_i$ for $i \in \braces{0, 1, \ldots, \omega}$ and all $\Gamma$ and $A : \Type_i(\Gamma)$.
We will then define a new cubical cwf $\MM^*$ denoted $\Con^*, \Subst^*, \ldots$, the \emph{Artin glueing} of $\MM$ with the standard model $\Std$ along an (internal) global sections functor, \ie the \emph{sconing} of $\MM$.
(We refrain from referring it to as just glueing to avoid confusion with the glue types of cubical cwfs.)

Recall from \cref{cubical-cwf} the operation $\ifill$ of $\MM$.
Instantiating it to the terminal context, we get $\Box \Fill(\Type(1), \lam{A} \Elem(1, A))$.
Using the forward direction of~\eqref{C-vs-fill-external}, we thus have an internal operation $k : \prod{A : \Type(1)} \CC(\Elem(1, A))$.

From now on, we will work in the internal language of presheaves over $\mathcal{C}$.
We start by defining a global sections operation $\verts{-}$ mapping contexts, types, and elements of $\MM$ to those of $\Std$.
\begin{itemize}
\item
Given $\Gamma : \Con$, we define $\verts{\Gamma} : \U_\omega$ as the set of substitutions $\Subst(1, \Gamma)$.
Given a substitution $\sigma : \Subst(\Delta, \Gamma)$, we define $\verts{\sigma} : \verts{\Delta} \to \verts{\Gamma}$ as $\verts{\sigma} \rho = \sigma \rho$.
This evidently defines a functor.
\item
Given $A : \Type(\Gamma)$, we define $\verts{A} : \verts{\Gamma} \to \Ufib_\omega$ as $\verts{A}\,\rho = (\Elem(1, A \rho), k\, (A\, \rho))$.
This evidently natural in $\Gamma$.
If $A$ is of level $n$, then $\verts{A} : \verts{\Gamma} \to \Ufib_n$.
\item
Given $a : \Elem(\Gamma, A)$ we define $\verts{a} : \prod{\rho : \Gamma}\,(\verts{A}\,\rho).1$ as $\verts{a}\,\rho = a \rho$.
This is evidently natural in $\Gamma$.
%This can be seen as the action of $\verts{-}$ on elements, but we choose not give a name to it.
\end{itemize}
Note that $\verts{-}$ preserves the terminal context and context extension up to canonical isomorphism in the category of contexts.
One could thus call $\verts{-}$ an (internal) \emph{pseudomorphism} cwfs from $\MM$ to $\Std$.
The sconing $\MM^*$ will be defined as essentially the Artin glueing along this pseudomorphism, but we will be as explicit as possible and not define Artin glueing at the level of generality of an abstract pseudomorphism.

For convenience, we also just write $\verts{A} : \verts{\Gamma} \to \U_\omega$ instead of $\lam{\rho} (\verts{A}\,\rho).1$, implicitly applying the first projection.
We also write just $\verts{A}$ for $\verts{A}\,\verts{()}$ if $\Gamma$ is the terminal context.

\subsection{Contexts, substitutions, types, and elements}
\label{sconing:contexts}

We start by defining the cwf $\MM^*$.

\begin{itemize}
\item
A context $(\Gamma, \Gamma') : \Con^*$ consists of a context $\Gamma : \Con$ in $\MM$ and a family $\Gamma'$ over $\verts{\Gamma}$ (which in the context of Artin glueing should be thought of as a substitution in $\Std$ from some context to $\verts{\Gamma}$).
We think of $\Gamma'$ as a \emph{proof-relevant computability predicate}.
A substitution $(\sigma, \sigma') : \Subst^*((\Delta, \Delta'), (\Gamma, \Gamma'))$ consists of a substitution $\sigma : \Delta \to \Gamma$ in $\MM$ and a map
\[
\sigma' : \prod{\nu : \verts{\Delta}} \Delta'(\nu) \to \Gamma'(\sigma \nu)
.\]
This evidently has the structure of a category.
\item
A type $(A, A') : \Type^*(\Gamma, \Gamma')$ consists of a type $A : \Type(\Gamma)$ in $\MM$ and
\[
A' : \prod{\rho : \verts{\Gamma}}\prod{\rho' : \Gamma'\,\rho} \verts{A}\,\rho \to \Ufib_\omega
.\]
We think of $A'$ as a fibrant \emph{proof-relevant computability family} on $A$.
In the abstract context of Artin glueing for cwfs, we should think of it as
\[
A' : \Type(\sum{\rho : \verts{\Gamma}}\sum{\rho' : \Gamma'\,\rho} \verts{A}\, \rho)
\]
in $\Std$.
However, we choose the former as the official definition so that the construction of $\MM^*$ from $\MM$ preserves Russell-style universes, as we shall see later.
Recalling $\Ufib_\omega = \sum{X : \U_\omega} C(X)$, we also write $\angles{A', \fib_{A'}}$ instead of $A'$ if we want to directly access the family and split off its proof of fibrancy.

The type $(A, A')$ is of level $n$ if $A$ and $A'$ are.

The action of a substitution $(\sigma, \sigma') : \Subst^*((\Delta, \Delta'), (\Gamma, \Gamma'))$ on $(A, A')$ is given by
\[
  (A \sigma, \lam{\nu,\nu',a} A'\,(\sigma \nu)\,(\sigma'\,\nu\,\nu')\, a).
\]
\item
An element $(a, a') : \Elem^*((\Gamma, \Gamma'), (A, \angles{A', \fib_{A'}}))$ consists of $a : \Elem(\Gamma, A)$ in $\MM$ and
\[
a' : \prod{\rho : \verts{\Gamma}}\prod{\rho' : \Gamma'\,\rho}\,A'(\rho, \rho', a \rho)
.\]
In the context of Artin glueing of cwfs (with types in $\MM^*$ presented correspondingly), this should be thought of as an element
\[
a' : \Elem(\sum{\rho : \verts{\Gamma}} \Gamma'\,\rho, \lam{(\rho, \rho')} A'(\rho, \rho', \verts{a}\,\rho))
\]
of $\Std$.

The action of a substitution $(\sigma, \sigma') : \Subst^*((\Delta, \Delta'), (\Gamma, \Gamma'))$ on the element $(a, a')$ is given by
\[
(a \sigma, \lam{\nu,\nu'} a'\,\sigma \nu\,(\sigma'\,\nu\,\nu'))
.\]
\item
The terminal context is given by $(1, 1')$ defined by $1'\,() = 1$.
\item
The extension in $\MM^*$  of a context $(\Gamma, \Gamma')$ by a type $(A, A')$ is given by $(\Gamma.A, (\Gamma.A)')$ where
\[
(\Gamma.A')(\rho, a) = \sum{\rho' : \Gamma'\,\rho}(A'\,\rho\,\rho'\,a).1
.\]
The projection $\pp^* : \Subst^*((\Gamma, \Gamma').(A, A'), (\Gamma, \Gamma'))$ is $(\pp, \pp')$ where
\[
\pp'\,(\rho, a)\,(\rho', a') = \rho'
\]
and the generic term $\qq^* : \Elem((\Gamma, \Gamma').(A, A') \pp^*)$ is $(\qq, \qq')$ where
\[
\qq'\,(\rho, a)\,(\rho', a') = a'
.\]
The extension of $(\sigma, \sigma') : \Subst^*((\Delta, \Delta'), (\Gamma, \Gamma'))$ with $(a, a') : \Elem^*((\Delta, \Delta'), (A, A')(\sigma, \sigma'))$ is
\[
((\sigma, a), \lam{\nu,\nu'} (\sigma'\,\nu\,\nu', a'\,\nu\,\nu'))
.\]
\end{itemize}

\subsection{Type formers and operations}

\subsubsection{Dependent products}

Let
\begin{align*}
(A, \angles{A', \fib_{A'}}) &: \Type^*(\Gamma, \Gamma')
,\\
(B, \angles{B', \fib_{B'}}) &: \Type^* ((\Gamma, \Gamma').(A, \angles{A', \fib_{A'}}))
.\end{align*}
We define the dependent product
\[
\iPi^*((A, \angles{A', \fib_{A'}}), (B, \angles{B', \fib_{B'}})) = (\iPi(A, B), \angles{\iPi(A, B)', \fib_{\iPi(A, B)'}})
\]
where
\[
\iPi(A, B)'(\rho, \rho', f) = \prod{a : \verts{A}\,\rho}\prod{a' : A'\,\rho\,\rho'\,a} B'\,(\rho, a)\,(\rho', a')\,(\iapp(f, a))
\]
and $\fib_{\iPi(A, B)'}(\rho, \rho', f)$ is given by closure of $\CC$ under dependent product applied to $(\verts{A}\,\rho).2$, $\fib_{A'}\,\rho\,\rho'\,a$ for $a : \verts{A}\,\rho$, and $\fib_{B'}\,(\rho, a)\,(\rho', a')\,(\iapp(f, a))$ for additionally $a' : A'\,\rho\,\rho'\,a$.

Given an element $(b, b')$ of $(B, \angles{B', d}))$ in $\MM^*$, we define the abstraction $\iabs^*(b, b') = (\iabs(b), \iabs(b)')$ where
\[
\iabs(b)'\,\rho\,\rho'\,a\,a' = b'\,(\rho, a)\,(\rho', a')
.\]

Given elements $(f, f')$ of $\iPi^*((A, \angles{A', c}), (B, \angles{B', d}))$ and $(a, a')$ of $(A, \angles{A', \fib_{A'}}))$ in $\MM^*$, we define the application $\iapp^*((f, f'), (a, a')) = (\iapp(f, a), \iapp(f, a)')$ where
\[
\iapp(f, a)'\,\rho\,\rho' = f'\,\rho\,\rho'\,a\rho\,(a'\,\rho\,\rho')
.\]

\subsubsection{Dependent sums}

Let
\begin{align*}
(A, \angles{A', \fib_{A'}}) &: \Type^*(\Gamma, \Gamma')
,\\
(B, \angles{B', \fib_{B'}}) &: \Type^* ((\Gamma, \Gamma').(A, \angles{A', \fib_{A'}}))
.\end{align*}
We define the dependent sum
\[
\iSigma^*((A, \angles{A', \fib_{A'}}), (B, \angles{B', \fib_{B'}})) = (\iSigma(A, B), \angles{\iSigma(A, B)', \fib_{\iSigma(A, B)'}})
\]
where
\[
\iSigma(A, B)'\,\rho\,\rho'\,(\ipair(a, b)) = \sum{a' : A'\,\rho\,\rho'\,a} B'\,(\rho, a)\,(\rho', a')\,b
\]
and $\fib_{\iSigma(A, B)'}\,\rho\,\rho'\,(\ipair(a, b))$ is given by closure of $\CC$ under dependent sum applied to $\fib_{A'}\,\rho\,\rho'\,a$ and $\fib_{B'}\,(\rho, a)\,(\rho', a')\,b$.

Given elements $(a, a')$ of $(A, \angles{A', \fib_{A'}})$ and $(b, b')$ of $(B, \angles{B', \fib_{B'}})[(a, a')]$ in $\MM^*$, we define the pairing $\ipair^*((a, a'), (b, b')) = (\ipair(a, b), \angles{a', b'})$.

Given an element $(\ipair(a, b), \angles{a', b'})$ of $\iSigma^*((A, \angles{A', \fib_{A'}}), (B, \angles{B', \fib_{B'}}))$ in $\MM^*$, we define the projections $\ifst^*(\ipair(a, b), \angles{a', b'}) = (a, a')$ and $\isnd^*(\ipair(a, b), {\angles{a', b'}}) = (b, b')$.

\subsubsection{Universes}

We define the universe $\iU_n^* : \Type^*(\Gamma, \Gamma')$ as $\iU_n^* = (\iU_n, \angles{\iU_n', \fib_{\iU_n'}})$ where
\[
\iU_n'\,\rho\,\rho'\,A = \verts{A}\,\rho \to \Ufib_n
\]
and $\fib_{\iU_n'}\,\rho\,\rho'\,A$ is given by $\CC(\Ufib_n)$ and closure of $\CC$ under exponentiation (note that fibrancy of $\verts{A}\,\rho$ is not used).
We have carefully chosen our definitions so that the evident natural isomorphism $\Elem^*((\Gamma, \Gamma'), \iU_n^*) \cong \Type_n^*(\Gamma, \Gamma')$ is an identity if the corresponding isomorphism in $\Type_n(\Gamma) \cong \Elem(\Gamma, \iU_n)$ in $\MM$ is an identity.
Thus, Russell-style universes are preserved by our presentation of the sconing model.

\subsubsection{Natural numbers}
\label{sconing:natural-numbers}

As per \cref{indexed-inductive-sets}, we have a fibrant indexed inductive set $\iN' \co \verts{\iN} \to \Ufib_0$ (where $\iN : \Type_0(1)$, hence $\verts{\iN} : \U_0$) with constructors
\begin{align*}
\zero' &: \iN'\,\izero,
\\
\succ' &: \prod{n : \verts{\iN}\,\rho} \iN'\,n \to \iN'\,(\isucc\,n)
.\end{align*}
%In the translation of \cref{indexed-inductive-sets}, the zero constructor unfolds for instance to $\zero'' : \Id_{\verts{\iN}}\,n\,\izero \to \iN'\,n$ for $n : \verts{\iN}$.
%This use of $\Id$ is critical to ensure fibrancy of $\iN'$.
In context $(\Gamma, \Gamma') : \Con^*$, we then define $\iN^* = (\iN, \lam{\rho,\rho'} \iN')$.
We have $\izero^* = (\izero, \lam{\rho,\rho'} \zero')$ and $\isucc^*(n, n') = (\isucc(n), \lam{\rho,\rho'} \succ'\,n \rho\,\ \!n')$ for $(n, n') : \Elem^*((\Gamma, \Gamma'), \iN^*)$.

Given $(P, P') : \Type((\Gamma, \Gamma').\iN^*)$ with
\begin{align*}
(z, z') &: \Elem^*((\Gamma, \Gamma')(P, P')[\izero^*])
,\\
(s, s') &: \Elem^*((\Gamma, \Gamma').\iN^*.(P, P'), (P, P') (\pp, S^*(\qq)) \pp)
\end{align*}
and $(n, n') : \Elem^*((\Gamma, \Gamma'), \iN^*)$, we define the elimination
\[
\inatrec^*((P, P'), (z, z'), (s, s'), (n, n')) = (\inatrec(P, z, s, n), \lam{\rho,\rho'} h'\,n \rho\,(n'\,\rho\,\rho'))
\]
where
\[
h' : \prod{m : \verts{\iN}}\prod{m' : \iN'\,m} P'\,(\rho, m)\,(\rho', m')\,(\inatrec(P \rho^+, z \rho, s \rho^{+++}, m))
\]
is given by induction on $\iN'$ with defining equations
\begin{align*}
h'\,\izero\,\zero' &= z'\,\rho\,\rho'
,\\
h'\,(\isucc(n))\,(\succ'\,n\,n') &= s'\,(\rho, n, \inatrec(P, z, s, n))\,(\rho', n', h'\,n\,n')
.\end{align*}

\subsubsection{Dependent paths}

Let $\angles{A, A'} : \II \to \Type^*(\Gamma, \Gamma')$ and $(a_b, a_b') : \Elem^*((\Gamma, \Gamma'), (A\,b, A'\,b))$ for $b \in \braces{0, 1}$.
We then define
\[
\iPath^*(\angles{A, A'}, (a_0, a_0'), (a_1, a_1')) = (\iPath(A, a_0, a_1), \angles{\iPath(A, a_0, a_1)', \fib_{\iPath(A, a_0, a_1)'}})
\]
where
\[
\iPath(A, a_0, a_1)'\,\rho\,\rho'\,(\ipabs(u)) = \Path_{\lam{i} (A'\,i\,\rho\,\rho'\,(u\,i)).1} (a'_0\,\rho\,\rho') (a'_1\,\rho\,\rho')
\]
and $\fib_{\iPath(A, a_0, a_1)'}\,\rho\,\rho'\,(\ipabs(u))$ is closure of $\CC$ under $\Path$ applied to $(A'\,i\,\rho\,\rho'\,(u\,i)).2$ for $i : \II$.

Given $\angles{u, u'} : \prod{i : \II}\Elem^*((\Gamma, \Gamma'), (A\,i, A'\,i))$, we define the path abstraction as
\[
\ipabs^*(\angles{u, u'}) = (\ipabs(u), \lam{\rho,\rho',i} u'\,i\,\rho\,\rho')
.\]

Given $(p, p') : \Elem^*((\Gamma, \Gamma'), \iPath^*(\angles{A, A'}, (a_0, a_0'), (a_1, a_1')))$ and $i : \II$, we define the path application
\[
\ipapp^*(p, i) = (\ipapp(p, i), \lam{\rho,\rho'} u'\,\rho\,\rho'\,i)
.\]

\subsubsection{Filling operation}

Given $\angles{A, A'} : \II \to \Type^*(\Gamma, \Gamma')$, $\phi : \FF$, $b \in \braces{0, 1}$, and
$\angles{u, u'} : \prod{i : \II}[\phi] \vee (i = b) \to \Elem^*((\Gamma, \Gamma'), (A\,i, A'\,i))$, we have to extend $u$ to
\[
\ifill^*(\angles{A, A'}, \phi, b, \angles{u, u'}) : \prod{i : \II} \Elem^*((\Gamma, \Gamma'), (A\,i, A'\,i))
.\]
We define $\ifill^*(\angles{A, A'}, \phi, b, \angles{u, u'}) = \angles{\ifill(A, \phi, b, u), \ifill(A, \phi, b, u)'}$ where
\[
\ifill(A, \phi, b, u)'\,i\,\,\rho\,\rho' : A'\,i\,\rho\,\rho'\,(\ifill(A, \phi, b, u)\,i)\rho
\]
is defined using $\fillop$ from~\eqref{U-fib-has-fill} as
\[
\ifill(A, \phi, b, u)'\,i\,\,\rho\,\rho' = \fillop(\lam{i} A'\,i\,\rho\,\rho'\,(\ifill(A, \phi, b, u)\,i)\rho, \phi, b, \lam{i,x} u'\,i\,x\,\rho\,\rho')
.\]

\subsubsection{Glue types}

Before defining the glueing operation in $\MM^*$, we will develop several lemmas relating notions such as contractibility and equivalences in $\MM$ with the corresponding notions of \cref{developments}.
Given $f : \Elem(\Gamma, A \to B)$ in $\MM$, we write
$
\efun{f} : \prod{\rho : \verts{\Gamma}} \verts{A}\,\rho \to \verts{B}\,\rho
$
for $\efun{f}\,\rho\,a = \iapp(f \rho, a)$.
This notation overlaps with the action of $\verts{-}$ on elements, but we will not use that one here.

Just in this subsection, we will use the alternative definition via given left and right homotopy inverses instead of contractible homotopy fibers of both equivalences $\iEquiv$ in the cubical cwf $\MM$ and equivalences $\Equiv$ in the (current) internal language.
In both settings, there are maps back and forth to the usual definition, which are furthermore natural in the context in the case of the cubical cwf $\MM$.
The statements we will prove are then also valid for the usual definition.

\begin{lem} \thmlabel{efun-pres-equiv}
  Given $f : \Elem(\Gamma, A \to B)$ in $\MM$ with $\Elem(\Gamma, \iisEquiv(f))$, we have
  $\prod{\rho : \verts{\Gamma}} \isEquiv(\efun{f}\,\rho)$.
This is natural in $\Gamma$.
\end{lem}

\begin{proof}
A (left or right) homotopy inverse $g : \Elem(\Gamma, B \to A)$ to $f$ in $\MM$ becomes a (left or right, respectively) homotopy inverse $\efun{g}\,\rho$ to $\efun{f}\,\rho$ for $\rho : \verts{\Gamma}$.
\end{proof}

\begin{lem} \thmlabel{equiv-in-sconing}
Given $(f, f') : \Elem((\Gamma, \Gamma'), (A, A') \to (B, B'))$ in $\MM^*$, the following statements are logically equivalent, naturally in $(\Gamma, \Gamma')$:
\begin{align}
\label{equiv-in-sconing:0}
&\Elem((\Gamma, \Gamma'), \iisEquiv^*(f, f'))
,\\
\label{equiv-in-sconing:1}
&\Elem(\Gamma, \iisEquiv(f)) \times \prod{\rho : \verts{\Gamma}}\prod{\rho' : \Gamma'\,\rho}
        \isEquiv(\sum{\efun{f}\,\rho} f'\,\rho\,\rho')
,\\
\label{equiv-in-sconing:2}
&\Elem(\Gamma, \iisEquiv(f)) \times \prod{\rho : \verts{\Gamma}}\prod{\rho' : \Gamma'\,\rho}\prod{a : \verts{A}\,\rho} \isEquiv(f'\,\rho\,\rho'\,a)
\end{align}
where $\sum{\efun{f}\,\rho} f'\,\rho\,\rho' : \sum{a : \verts{A}\,\rho} A'\,\rho\,\rho'\,a \to
\sum{b : \verts{B}\,\rho} B'\,\rho\,\rho'\,b$.
\end{lem}

\begin{proof}
Let us only look at homotopy left inverses.

For~$\eqref{equiv-in-sconing:0} \to \eqref{equiv-in-sconing:1}$, a homotopy left inverse $(g, g')$ to $(f, f')$ in $\MM^*$ gives a homotopy left inverse $\sum{\efun{g}\,\rho} g'\,\rho\,\rho'$ to $\sum{\efun{f}\,\rho} f'\,\rho\,\rho'$ for all $\rho, \rho'$.

For~$\eqref{equiv-in-sconing:1} \to \eqref{equiv-in-sconing:2}$, we use \cref{efun-pres-equiv} and note that a fiberwise map over an equivalence is a fiberwise equivalence exactly if it is an equivalence on total spaces (the corresponding statement for identity types instead of paths is~\cite[Theorem~4.7.7]{UF}).

For~$\eqref{equiv-in-sconing:2} \to \eqref{equiv-in-sconing:0}$, given a homotopy left inverse $g$ to the equivalence $f$ in $\MM$ and a homotopy left inverse $\overline{g}'\,\rho\,\rho'\,a : B'\,\rho\,\rho'\,(\efun{f}\,a) \to A'\,\rho\,\rho'\,a$ to $f'\,\rho\,\rho'\,a$ for all $\rho, \rho', a$, we use \cref{efun-pres-equiv} to transpose $\overline{g}'$ to the second component $g'\,\rho\,\rho'\,b : B'\,\rho\,\rho'\,b \to A'\,\rho\,\rho'\,(\efun{g}\,b)$ for all $\rho, \rho', b$ of a homotopy left inverse $(g, g')$ to $(f, f')$ in $\MM^*$.
\end{proof}

We can now define glue types in $\MM^*$.
Let $(A, A') : \Type(\Gamma, \Gamma')$, $\phi : \FF$, $\angles{T, T'} : [\phi] \to \Type(\Gamma, \Gamma')$, and
\[
\angles{e, e'} : [\phi] \to \Elem((\Gamma, \Gamma'), \iEquiv^*((T \TT, T' \TT), (A, A')))
.\]
We define
\[
\iGlue^*((A, A'), \phi, \angles{T, T'}, \angles{e, e'}) = (\iGlue(A, \phi,, T, e), \angles{G', \fib_{G'}})
\]
where
\[
G'\,\rho\,\rho'\,(\iglue(a, t)) = \Glue\,(A'\,\rho\,\rho'\,a).1\,[\phi \mapsto (T'\,\TT\,\rho\,\rho'\,(t\,\TT), ((e'\,\TT\,\rho\,\rho').1\, (t \,\TT)))]
\]
where $\fib_{G'}\,\rho\,\rho'\,(\iglue(a, t))$ is given by closure of $\CC$ under $\Glue$ applied to $(A'\,\rho\,\rho'\,a).2$ and $T'\,\TT\,\rho\,\rho'\,(t\,\TT)$ on $[\phi]$ and the witness that $(e'\,\TT\,\rho\,\rho').1 \,(t \,\TT)$ is an equivalence provided by the direction from~\eqref{equiv-in-sconing:0} to~\eqref{equiv-in-sconing:2} of \cref{equiv-in-sconing}.
We define $\iunglue^* = (\iunglue, \iunglue')$ where
\[
\iunglue'\,\rho\,\rho'\,(\iglue(a, t)) = \unglue
.\]
Given $(a, a') : \Elem((\Gamma. \Gamma'), (A, A'))$ and
\[
(t, t') : [\phi] \to \Elem((\Gamma, \Gamma'), (T\,\TT, T'\,\TT))
\]
such that $\iapp^*(\ifst^*(e, e')\,\TT, (t, t')\,\TT) = (a, a')$ on $[\phi]$, we define $\iglue^*((a, a'), (t, t'))$ as the pair $(\iglue(a, t), \iglue(a, t)')$ where
\[
\iglue(a, t)'\,\rho\,\rho' = \glue\,(a'\,\rho\,\rho')\,[\phi \mapsto t'\,\TT\,\rho\,\rho']
.\]

\subsection{Main result}

One checks in a mechanical fashion that the operations we have defined above satisfy the required laws, including stability under substitution in the context $(\Gamma, \Gamma')$.
We thus obtain the following statement.

\begin{thm}[Sconing] \thmlabel{sconing-result}
Assume the parameters $\mathcal{C}, \II, \FF$ satisfy the assumptions of \cref{developments}.
Then given any cubical cwf $\MM$ that is size-compatible in the sense of the beginning of \cref{sconing}, the sconing $\MM^*$ is a cubical cwf with operations defined as above.
We further have a morphism $\MM^* \to \MM$ of cubical cwfs given by the first projection.
\end{thm}

\section{Homotopy canonicity}

We fix parameters $\mathcal{C}, \II, \FF$ as before.
To make our homotopy canonicity result independent of \cref{thm:initial-ccwf} concerning initiality of the term model, we phrase it directly using the \emph{initial model} $\mathcal{I}$, initial in the category of cubical cwfs with respect to the parameters $\mathcal{C}, \II, \FF$.
Its existence can be justified generically following~\cite{sterling:algebraic-type-theory,palmgren-vickers:partial-horn-logic}.
It is size-compatible in the sense of \cref{sconing}: internally, $\Subst_{\mathcal{I}}(\Delta, \Gamma)$ and $\Elem_{\mathcal{I}}(\Gamma, A)$ live in the lowest universe $\U_0$ for all $\Gamma, \Delta, A$.

\begin{thm}[Homotopy canonicity]
Assume the parameters $\mathcal{C}, \II, \FF$ satisfy the assumptions of \cref{developments}.
In the internal language of presheaves over $\mathcal{C}$, given a closed natural $n : \Elem(1, \iN)$ in the
initial model $\mathcal{I}$, we have a numeral $k : \N$ with $p : \Elem(1, \iPath(\iN, n, \isucc^k(\izero)))$.
\end{thm}

\begin{proof}
We start the arguing reasoning externally.
Using \cref{sconing-result}, we build the sconing $\mathcal{I}^*$ of $\mathcal{I}$.
Using initiality, we obtain a section $F$ of the cubical cwf morphism $\mathcal{I}^* \to \mathcal{I}$.

Let us now proceed in the internal language.
Recall the construction of \cref{sconing:natural-numbers} of natural numbers in $\mathcal{I}^*$.
We observe that $\sum{n : \verts{\iN}} \iN'\,n$ forms a fibrant natural number set (in the sense of \cref{indexed-inductive-sets}).
It is thus homotopy equivalent to $\N$.
Under this equivalence, the first projection $\sum{n : \verts{\iN}} \iN'\,n \to \verts{\iN}$ implements the map sending $k : \N$ to $\isucc^k(\izero).$

Inspecting the action of $F$ on $n : \Elem(1, \iN)$, we obtain $n' : \iN'\,n$.
By the preceding paragraph, this corresponds to $k : \N$ with a path $p' : \II \to \verts{\iN}$ from $n$ to $\isucc^k(\izero)$.
Now $p = \ipabs(p')$ is the desired witness of homotopy canonicity.
\end{proof}

%% \subsection{A concrete example}

%% Since the model is effective, we can produce, given $t$ of type $N_2$,
%% an {\em actual} value $0$ or $1$ and
%% a path in $|N_2|$ between $t$ and this value.

%%  A concrete example is the following. We can consider the equivalence $\neg:N_2\to N_2$
%% defined by the negation.
%% We can then consider the type $T(i) = \iGlue~[i=0\mapsto (N_2,id),~i=1\mapsto (N_2,\neg)]~N_2$
%% which is a path between $T(0) = N_2$ and $T(1) = N_2$
%% and the term $t = \comp^i~T(i)~[]~0$ of type $N_2$. We compute an element in ${N_2}'(t)$ and
%% this produces the element $1$ and a path between $1$ and $t$ in $|N_2|$.
%% This computation is possible since ${N_2}'$ defines a fibration over $|N_2|$.

\section{Extensions}

\subsection{Identity types}
\label{identity-types}

Our treatment extends to the variation of cubical cwfs that includes identity types.

\emph{Identity types} in a cubical cwf denoted as in \cref{cubical-cwf} consist of the following operations and laws (omitting stability under substitution), internal to presheaves over $\mathcal{C}$.
Fix $A$ in $\Type(\Gamma)$.
Given $x, y$ in $\Elem(\Gamma, A)$, we have $\iId(A, x, y)$ in $\Type(\Gamma)$, of level $n$ if $A$ is.
Given $a$ in $\Elem(\Gamma, A)$, we have $\irefl(a)$ in $\Elem(\Gamma, \iId(A, a, a))$.
Given $P$ in $\Type(\Gamma.A.A \pp.\iId(A \pp \pp, \qq \pp, \qq))$ and $d$ in $\Elem(\Gamma.A, P[\qq, \qq, \irefl(q)])$ and $x, y$ in $\Elem(\Gamma, A)$ and $p$ in $\Elem(\Gamma, \iId(A, x, y))$, we have $\iJ(P, d, x, y, p)$ in $\Elem(\Gamma, P[x, y, p])$.
We have
\[
\iJ(P, d, a, a, \irefl(a)) = d[a]
.\]

We can interpret univalent type theory in any cubical cwf with identity types as per \cref{univalent-type-theory-prop-J}.

The standard model of \cref{standard-model} has identity type $\iId(\angles{A, \fib_A}, x, y) : \Type(\Gamma)$ given by
$\prod{\rho : \Gamma} \Id_{A\,\rho}\,(x\,\rho)\,(y\,\rho)$ using Andrew Swan's construction of $\Id$ referenced in \cref{indexed-inductive-sets}.
We omit the evident description of the remaining operations.

To obtain homotopy canonicity in this setting, it suffices to extend the sconing construction $\MM^*$ of \cref{sconing} to identity types.
Given $A : \Type(1)$ and $A' : \verts{A} \to \Ufib_\omega$, we define $\Id'_{A, A'}$ as the fibrant indexed inductive set (as per \cref{indexed-inductive-sets}) over $x, y : \verts{A}$, $p : \verts{\iId(A, x, y)}$, $x' : A'\,x$, $y' : A'\,y$ with constructor
\[
\refl' : \prod{a : \verts{A}}\prod{a' : A'\,a} \Id'_{A, A'}\,a\,a\,(\irefl(a))\,a'\,a'
.\]

Now fix $(A, A') : \Type^*(\Gamma, \Gamma')$.
Given $\rho : \verts{\Gamma}$, $\rho' : \Gamma'\,\rho$, and
elements $(x, x'), (y, y')$ of $(A, A')$ in $\MM^*$, we define
\[
\iId^*((A, A'), (x, x'), (y, y')) = (\iId(A, x, y), \lam{\rho,\rho',p} \Id'_{A \rho, A'\,\rho\,\rho'}\,x \rho\ y \rho\ p\,(x'\,\rho\,\rho')\,(y'\,\rho\,\rho'))
.\]
Given an element $(a, a')$ of $(A, A')$ in $\MM^*$, we define $\irefl^*(a, a') = (\irefl(a), \irefl(a)')$ where
\[
\irefl(a)'\,\rho\,\rho' = \refl'\,a \rho\ (a'\,\rho\,\rho')
.\]
The eliminator $\iJ((C, C'), (d, d'), (x, x'), (y, y'), (p, p'))$ is defined as
\[
(\iJ(C, d, x, y, p), \lam{\rho,\rho'} h'\,x \rho\ y \rho\ p \rho\ (x'\,\rho\,\rho')\,(y'\,\rho\,\rho')\,(p'\,\rho\,\rho'))
\]
where
\[
\begin{split}
h' : &\prod{x, y : \verts{A}\,\rho} \prod{p : \verts{\iId(A, x, y)}\,\rho} \prod{x' : A'\,\rho'} \prod{y' : A'\,\rho'} \prod{p' : \Id'_{A \rho, A'\,\rho\,\rho'}\,x\,y\,p\,x'\,y'}
\\
& P'\,(\rho, x, y, p)\,(\rho', x', y, p')\,(\iJ(P \rho^{+++}, d \rho^+, x, y, p))
\end{split}
\]
is given by induction on $\Id'_{A \rho, A'\,\rho\,\rho'}$ via the clause
\[
h'\,a\,a\,(\irefl(a))\,a'\,a'\,(\refl'\,a\,a') = d'\,(\rho, a)\,(\rho', a')
.\]

\subsection{Higher inductive types}
\label{higher-inductive-types}

Our treatment extends to higher inductive types~\cite{UF}, following the semantics presented in~\cite{CHM}.
Crucially, we have fibrant \emph{indexed} higher inductive sets in presheaves over $\mathcal{C}$ as we have what we would call fibrant \emph{uniformly indexed} higher inductive sets in the same fashion as in~\cite{CHM} and fibrant identity sets~\cite{CCHM,OP}, mirroring the derivation of fibrant indexed inductive sets from fibrant uniformly indexed inductive sets and fibrant identity sets recollected in \cref{indexed-inductive-sets}.%
\footnote{We stress that the use of ``set'' in this context refers to the types of the language of presheaves over $\mathcal{C}$, not homotopy sets.} %

Let us look at the case of the \emph{suspension} operation in a cubical cwf, where $\isusp(A) : \Type(\Gamma)$ has constructors $\inorth, \isouth$ and $\imerid(a, i)$ for $a : A$ and $i : \II$ with $\imerid(a, 0) = \inorth$ and $\imerid(a, 1) = \isouth$.

For the sconing model of \cref{sconing}, we define for $A : \Type(1)$ and $A' : \verts{A} \to \Ufib_\omega$ the indexed higher inductive set $\susp'_{A, A'}$ over $\verts{\isusp(A)}$ with constructors
\begin{align*}
\north' &: \susp'_{A, A'}\,\inorth
,\\
\south' &: \susp'_{A, A'}\,\isouth
,\\
\merid'\,a\,a'\,i &: (\susp A)'(\imerid(a, i))[i=0\mapsto \north',i=1\mapsto \south']
\end{align*}
for $a : \verts{a}$ and $a' : A'\,a$ and $i : \II$ (using the notation of~\cite{CHM}).
In the above translation to a uniformly indexed higher inductive set, the constructor $\north'$ will for example be replaced by
\[
\north'' : \Id_{\verts{\isusp(A)}}\,u\,\inorth \to \isusp'_{A, A'}\,u
.\]

Given $(A, A') : \Type^*(\Gamma, \Gamma')$, we then define
\[
\isusp^*(A, A') = (\isusp(A), \lam{\rho,\rho'} \susp'_{A\rho, A'\,\rho\,\rho'})
,\]
with constructors and eliminator treated as in \cref{identity-types}.

\section{Canonicity}
\label{strict-canonicity}

The goal of this section is to show \emph{canonicity} for cubical type theory, stating that any closed term of type $\iN$ is (strictly) equal to a numeral.
This is a priori a stronger result than merely \emph{homotopy canonicity}.
However, it requires us to add further computation rules for the filling operation to the theory.
For this purpose, we have defined in \cref{computational-cubical-cwf} the notion of computational cubical cwf, modelling a modified version of cubical type theory where the filling operation is replaced by the composition operation (filling is then a {\em derived} operation).
This is our notion of model in this section.

The main point is how to define the right notion of computability structure.
Once this is done, we can essentially construct the sconing model as in~\cite{COQ2019}.
If we apply this to the initial computational cubical cwf, we get the canonicity
result: any closed natural number term is convertible to a numeral.
This result was already proved in~\cite{Huber}, but
like for the proof in~\cite{COQ2019}, our new argument completely avoids the
need to define a reduction relation, which is quite subtle for cubical
type theory in~\cite{Huber} since it is not closed under name substitution.

As in \cref{computational-cubical-cwf}, we have not just connection structure on the interval, but also a compatible reversal structure.
Other than that, we make the same assumptions on the $\II$ and $\FF$ as in \cref{developments}.
Starting from an arbitrary computational cubical cwf $\MM$ in the global context satisfying size-compatibility as in \cref{sconing}, we build a new computational cwf $\MM^*$, the sconing of $\MM$.

\subsection{Sconing model: cwf structure}
\label{canonicity-sconing-cwf}

The underlying category of $\MM^*$ is defined as in \cref{sconing}.
It is the Artin glueing of $\MM$ along the global sections functor $\verts{-} \colon \MM \to \U_\omega$.
(Recall that $\verts{\Gamma}$ is $\Subst(1, \Gamma)$ for $\Gamma$ in $\Con$.
In particular, a context in $\MM^*$ is a pair $(\Gamma, \Gamma')$ where $\Gamma$ is a context in $\MM$ and $\Gamma' : \verts{\Gamma} \to \U_\omega$.)

In contrast to \cref{sconing}, we cannot view the rest of the structure $\MM^*$ as being obtained by glueing along the pseudomorphism $\verts{-} \colon \MM \to \Std$ of (computational) cubical cwfs.
In particular, we will not make use the standard model $\Std$.
Rather, $\MM^*$ can be seen as the total space of a fibration that presents a fibred version of the standard model over $\MM$.
In the definition of types, we need to track computability of the composition operation, expressed using the right adjoint to exponentiation with $\II$ to get a fiberwise notion.

For $A$ in $\Type(1)$, we write $\verts{A}$ for $\Elem(1, A)$.
Set $\Pred_n = \sum{A:\Type(1)} \verts{A} \to \U_n$.
For $\lam{i} (A_i,A'_i)$ in $\Pred_{\omega}^\II$, we define $\Red(\lam{i} (A_i,A'_i))$ to be the type of operations $c$ taking as argument $\psi$ in $\FF$ and a family $u_i,u'_i$ for $[\psi] \vee (i = 0)$ with $u_i$ in $A_i$ and $u'_i$ in $A'_i\,u_i$ and producing an element $c(\psi,u,u')$ in $A'_1(\icomp(A,\psi,u))$
which is equal to $u'_1$ for $\psi = 1$.

\begin{prop} \thmlabel{defRed}
We have $R:\Pred_\omega\to\U_\omega$ with, naturally in global $A$ and $B:A\to\Pred_{\omega}$, a bijection between $\Box(\prod{A} R \circ B)$ and $\Box(\prod{A^{\II}} \Red \circ B^\II)$.
For $n \geq 0$, we have that $R$ descends to an operation $R : \Pred_n \to \U_n$.
\end{prop}

\begin{proof}
This follows from \cref{adjunction-slice}.
\end{proof}

\begin{prop} \thmlabel{mainth}
Given $e$ in $\prod{i:\II} R(A_i,A'_i)$, there is an operation $c(e)$ which, given $\psi$ in $\FF$ and $u'_i$ in $A'_i\,u_i$ for $[\psi] \vee (i = 0)$, produces an element $c(e)(\psi,u,u')$ in $A'_1(\icomp(\lam{i} A_i,\psi,u))$ which is equal to $u'_1$ for $[\psi]$.
\end{prop}

\begin{proof}
This follows from \cref{defRed} by setting $T = \sum{\Pred_\omega}R$.
\end{proof}

\begin{rem}
Define $C(X)$ for $X : \U_\omega$ as in \cref{fibrant}, but using composition instead of filling.
As in \cref{sconing}, we have a map $k : \prod{A : \Type(1)} C(|A|)$.
In the same fashion as above for $R$, one may construct $C'$ over $\sum{A : \U_\omega} \sum{c : C(A)} {A' : A \to \U_\omega}$ encoding that the family $A'$ has ``composition over $c$''.
Then $R$ can be defined as the restriction of $C'$ along the map induced by $k$.
Under the equivalence between families $A' : A \to \U_\omega$ over $A$ and $\overline{A} : \U_\omega$ with a map $p : \overline{A} \to A$, this corresponds to an element of $C(\overline{A})$ such that $p$ forms a composition-preserving morphism of fibrant types, a notion defined (like $C$ and $C'$) using the adjunction recorded in \cref{adjunction-slice}.

The closure properties of $R$ under type formers proved in \cref{canonicity-dependent-sum,nat-computable} below can be proved at the level of $C'$.
This has the advantage of eliminating the dependency on the computational cubical cwf $\MM$ (in the case of natural numbers, the input is instead a natural number algebra in fibrant types).
Then the actual sconing construction can proceed without external reasoning as in \cref{sconing}.

The sconing of the standard model for computational cwfs (defined as $\Std$ in \cref{standard-model}, but using composition instead of filling) will coincide, externally, with the standard model constructed internally in presheaves over $[1] \times \mathcal{C}$ (where $[1]$ is the poset with elements $0 < 1$) with interval object and cofibration classifier defined by projection to $\mathcal{C}$.
\end{rem}

We define $\Type^*(\Gamma,\Gamma')$ to be the set of triples $(A,A',e_A)$ where $A$ is in $\Type(\Gamma)$ and $A'\,\rho\,\rho'$ is in $\verts{A\,\rho} \to \U_\omega$ for $\rho$ in $\verts{\Gamma}$ and $\rho'$ in $\Gamma'\,\rho$ and $e_A\,\rho\,\rho'$ is in $R(A\rho,A\,\rho\,\rho')$ for $\rho$ and $\rho'$ as before.
%\note{Define an auxiliary notion of semantic type over $\U_\fib$ that groups $A'$ and $e_A$?}
%\note{The name $e_A$ is a bit bad because it depends on $A'$.}

We define $\Elem^*((\Gamma,\Gamma'),(A,A',e_A))$ to be the set
of pairs $(a,a')$ where $a$ is in $\Elem(\Gamma,A)$ and
$a'\,\rho\,\rho'$ is in $A'\,\rho\,\rho'\,(a\rho)$ for $\rho$ in $\verts{\Gamma}$ and
$\rho'$ in $\Gamma'\,\rho$.

Since the elements do not make use of the last component of the triple defining a type, the operations involving context extension are defined as in \cref{sconing:contexts}.
For example, we define $(\Gamma,\Gamma').(A,A',e_A)$ to be $(\Gamma.A,(\Gamma.A'))$ where $(\Gamma.A)'(\rho,u)$ is $\sum{\rho':\Gamma'\,\rho} A'\,\rho\,\rho'\,u$.

%\cref{mainth} can then be used to give an interpretation of the composition operation in this model.

\subsection{Example of dependent sum types}
\label{canonicity-dependent-sum}

Before explaining the example of the dependent sum type, we need the following preliminary lemma.
It intuitively says that the filling operation is computable if the composition operation is computable.

\begin{lem} \thmlabel{compFill}
Given $e$ in $\prod{i:\II} R(A_i,A'_i)$, the filling operation on $A$ is ``computable'': for $\psi$ in $\FF$ and a partial family $u_i$ of elements in $\verts{A_i}$ together with $u_i'$ in $A'_i\,u_i$ defined for $[\psi] \vee (i=0)$, then for any $r$ in $\II$ we have
\[
\ifill'(A,\psi,u,u')\,r : A'_r\,(\ifill(A,\psi,u)\,r)
\]
equal to $u'_r$ for $[\psi] \vee (r = 0)$.
\end{lem}

%\note{Was coherence of $\ifill$ with $c(e)$ needed later?}

\begin{proof}
Given $r$ in $\II$, we define $e_r$ in $\prod{i:\II} R(A_{i\wedge r},A'_{i\wedge r})$ by $e_r\,i = e_{r\wedge i}$.
Using \cref{mainth}, we can define
\[
\ifill'(A,\psi,u,u')\,r = c(e_r)(\psi,u_r,u_r')
\]
with $u_r\,i\,x = u\,(i \wedge r)\,x$ and $u_r'\,i\,x = u'\,(i \wedge r)\,x$.
\end{proof}

Given $(A,A')$ in $\Pred_\omega$ and $B:\Type(1.A)$ and $B' : \prod{u:\verts{A}} A'\,u \to \verts{B[u]} \to \U_\omega$
we define $\iSigma(A,B)'$ by
\[
\iSigma(A,B)'\,w = \sum{u':A'\,(\ifst(w))}\,B'\,(\ifst(w))\,u'\,(\isnd(w))
.\]
We then want to define an operation
\[
R(A,A') \to (\prod{u:\verts{A},u':A'\,u} R(B[u],B'\,u\,u')) \to R(\iSigma(A,B),\iSigma(A,B)')
.\]
In order to do this, we consider the iterated dependent sum $\Delta$
corresponding to the context
\begin{align*}
(A,A') &: \Pred_n,
\\
B &: \Type_n(A),
\\
B' &: \prod{u:\verts{A}} A'\,u \to \verts{B[u]} \to \U_n,
\\
e_A &: R(A,A'),
\\
e_B &: \prod{u:\verts{A}} \prod{u':A'\,u} R(B[u],B'\,u\,u')
.\end{align*}
We want to build a global element of
\[
\prod{(A,A',B,B',e_A,e_B):\Delta}\,R(\iSigma(A,B),\iSigma(A,B)')
.\]
Using \cref{defRed}, we are reduced to show the following statement.

\begin{prop}
There is an element of
\[
\prod{\lam{i} (A_i,A'_i,B_i,B'_i,e_A^i,e_B^i):\Delta^\II}
\Red(\lam{i} (\iSigma(A_i,B_i),\iSigma(A_i,B_i)'))
.\]
\end{prop}

\begin{proof}
We assume a family $A_i,A'_i$ in $\Pred_n$ and $e_A^i$ in $R(A_i,A'_i)$ for $i : \II$.
We also have $B_i:\Type_n(A_i)$ and $B'_i$ in $\prod{u:\verts{A_i}} A_i'\,u \to \verts{B_i[u]} \to \U_n$ and $e_B^i\,u\,u'$ in $R(B_i[u],B'_i\,u\,u')$.
We also have $\psi$ in $\FF$ with $\ipair(a_i,b_i) : \iSigma(A_i,B_i)$ and $a'_i$ in $A'_i\,a_i$ and $b'_i$ in $B'_i\,a_i\,a'_i\,b_i$ defined for $[\psi] \vee (i = 0)$.

Given all this, we want to build
\[
(a'_1,b'_1) : \iSigma(A_1,B_1)'(\icomp(\iSigma(A,B),\psi,(a,b)))
.\]
By the computation rule for $\iSigma$, we have
\[
\icomp(\iSigma(A,B),\psi,(a,b)) = (\icomp(A,\psi,a),\icomp(\lam{i} B[u\,i],\psi,b))
\]
where $u = \ifill(A,\psi,a)$.

Using \cref{compFill}, we have $u' = \ifill'(A,\psi,a,a')$ and we define $a'_1 = u'\,1$.
We next define $e_B\,i = e_B^i\,(u\,i)\,(u'\,i)$, and using \cref{mainth}, we take $b'_1 = c(e_B)(\psi,b,b')$.
\end{proof}

We can use these results to interpret dependent sum types in the model $\MM^*$, following
essentially the interpretation in~\cite{COQ2019}.

Given $(A,A',e_A)$ in $\Type^*(\Gamma,\Gamma')$ and
$(B,B',e_B)$ in $\Type^*((\Gamma,\Gamma').(A,A',e_A))$
we define $T,T',e_T$ where $T = \iSigma(A,B)$ in $\Type^*(\Gamma,\Gamma')$.

Given $\rho$ and $\rho'$ we have $A_1 = A\rho$ in $\Type(1)$ and
$A_1' = A'\,\rho\,\rho'$ in $\verts{A_1} \to \U_\omega$. We also have $B_1 = B\rho^+$ in
$\Type(A_1)$ and $B_1 u u' = B'(\rho,\ifst(w))(\rho',u')$
in
$\prod{u:\verts{A_1}}A_1'u\to \verts{B_1[u]} \to \U_\omega$
since $B_1[u] = B\rho^+[u] = B(\rho,u)$. We can then define $T'\,\rho\,\rho'$
to be $\iSigma(A_1,B_1)'$ and use $(1)$ to define $e_T\,\rho\,\rho'$.

\subsection{Natural numbers}
\label{canonicity-nat}

Let $\NN$ be the (internal) set of natural numbers (given by the constant presheaf of natural numbers).
We have a canonical map $\QUOTE:\NN \to \verts{\iN}$ sending $k$ to $\isucc^k(\izero)$.
We define a (non-fibrant) family $\iN'$ over $\verts{\iN}$ by $\iN'(t) = \sum{k:\NN}t =_{\verts{\iN}}\QUOTE(k)$ (using the {\em strict} equality on the set $\verts{\iN}$).%
\footnote{An isomorphic alternative is to define $\iN'$ as a (\emph{non-fibrant}) indexed inductive set in the presheaf model, with constructors of type $\iN'\,\izero$ and $\iN'\,n \to \iN'\,(\isucc\,n)$ for $n : \verts{\iN}$.
Indeed, it is this approach that generalizes to the interpretation of inductive types with parameters.}

\begin{lem} \thmlabel{nat-computable}
We have an element of $R(\iN,\iN')$.
\end{lem}

\begin{proof}
Using the adjoint definition of $R$ in \cref{defRed}, we must build $c : \Red(\lam{i} (\iN,\iN'))$.
Exponentiation with $\II$ preserves external coproducts since $\II$ is tiny.
Since $\NN$ is a countable coproduct of $1$, it follows that any function $\II \to \NN$ is constant (formally, factors uniquely through $\II \to 1$).

Let $u_i : |\iN|$ and $u'_i : \iN'\,u_i$ for $[\psi] \vee (i = 0)$.
The latter means $k_i : \NN$ such that $u_i = \QUOTE(k_i)$ for $[\psi] \vee (i = 0)$.
Using the observation from the previous paragraph, there is unique $k : \NN$ such that $k_i = k$ for $[\psi] \vee (i = 0)$.
From the equations for the composition operation on $\iN$ in a computational cubical cwf and induction on $k$, we get that
\[
\icomp(\lam{i} \iN, \psi, \lam{i, x} u_i) = \icomp(\lam{i} \iN, \psi, \lam{i, x} \QUOTE(k)) = \QUOTE(k)
.\]
We are forced to set $c(\psi, u, u') = (k, \TT)$.
\end{proof}

This provides the interpretation of the type of natural numbers in the model $\MM^*$.

We see here a key difference compared to the sconing model used for proving homotopy canonicity: the computability predicate used in this case is not valued in fibrant sets.
Note that the family $\sum{k:\NN} \verts{\Path(\iN,t,\QUOTE(k))}$ for $t:\verts{\iN}$ is fibrant, but would not work for showing canonicity since it does not support an interpretation of $\inatrec$.

\subsection{Proof of canonicity}

Starting from any computational cubical cwf $\MM$, we have
built a new model $\MM^*$, the associated {\em computability model},
with a (strict) projection map $\MM^*\to\MM$. Like in~\cite{COQ2019},
if we apply this to the initial model, we get that any closed term of type $\iN$
is ``computable'', \ie, is {\em strictly} equal to a numeral.

\section*{Conclusion}

We have given proofs of two forms of canonicity for cubical type theory.
The first one is homotopy canonicity (every closed term of type $\iN$ is {\em path equal} to a numeral) in a cubical type theory without structural computation rules for the composition operation.
The second one is canonicity (every closed term of type $\iN$ is \emph{strictly equal} to a numeral) in a cubical type theory with these computation rules.
While our arguments rely on an interplay between internal and external reasoning, the main part of the first argument can be seen as happening internally in the model of fibrant sets.
The second argument can hopefully be refined to a constructive proof of normalisation.

\bibliographystyle{alphaurl}
\bibliography{references}

\appendix

\section{Rules of the term model}
\label{term-model-rules}

We denote the objects of our base category $\mathcal{C}$ by $X,Y,Z$
and its morphisms by $f,g,h$.  In the term model $\Syn$ morphisms $f
\co Y \to X$ act on judgments at stage $X$ via an implicit
substitution, while for substitutions on object variables we will use
explicit substitutions.  For this to make sense we first define the
raw expressions as a presheaf: at stage $X$ this is given by
\[
  \begin{array}{rcl}
    \Gamma,\Delta &::=
    & \emptyctx \mid \Gamma.A
    \\
    A,B,t,u,v & ::=
    & \qq \mid t \sigma \mid
    % univ
      \iU_n \mid
      % pi
      \iPi (A,B) \mid \iabs (u) \mid \iapp (u,v)
    \\
                  & \mid
    &
    % sigma
      \iSigma (A,B) \mid \ipair (u,v) \mid \ifst(u)
      \mid \isnd (u)
      % path
    \\
                  & \mid
    &
      \iPath (\bar A,u,v) \mid \ipabs \bar u \mid \ipapp (u,r)
    \\
                  & \mid
    & \iGlue (A, \phi, \bar B, \bar u)
      \mid \iglue (v,\bar u) \mid \iunglue (u)
    \\
                  & \mid
    & \ifill (\bar A, \phi, b, \bar u, r)
      \mid \dots
    \\
    \bar A, \bar B, \bar u, \bar v
                  &::=& (A_{f,r})_{f,r} \mid (A_f)_{f \in [\phi]}
    \\
    \sigma,\tau,\delta
                  &::=& \pp \mid \id \mid \sigma \tau \mid
                        (\sigma,u) \mid \emptysubst
    % \\
    % r & \in & \II (X)
    % \\
    % \phi & \in & \FF (X)
  \end{array}
\]
where $b \in \braces{0,1}$, $\phi \in \FF (X)$, and we skipped the
constants for natural numbers.  Above, we have families of
expressions, say $\bar A = (A_{f,r})_{f,r}$, whose index set ranges
over certain $Y$, $f \co Y \to X$, and $r \in \II(Y)$, and $A_{f,r}$
is a raw expression at stage $Y$; likewise $(A_f)_{f \in [\phi](X)}$
consists of raw expressions $A_f$ at stage $Y$ for $f \co Y \to X$ in
the sieve $[\phi]$ on $X$.  (The exact index sets will be clear from
the typing rules below.)  All other occurrences of $r$ above have $r
\in \II(X)$.  The restrictions along $f \co Y \to X$ on the raw syntax
then leave all the usual cwf structure untouched, so we have $\qq f =
\qq$ and $(\iPi (A,B))f = \iPi (A f, Bf)$, and uses the restrictions
in $\II$ and $\FF$ accordingly, \eg, $(\ipapp (u,r)) f = \ipapp (uf ,
r f)$, and we will re-index families according to ${\bar A} f =
(A_{gf,rf})$ for $\bar A = (A_{g,r})_{g,r}$.

To get the \emph{initial} cubical cwf we in fact need more annotations
to the syntax in order to be able to define a partial interpretation
(\cf~\cite{Streicher91,Hofmann97}) on the raw syntax.  But to enhance
readability we suppress these annotations.

We will now describe a type system indexed by stages $X$.  The forms
of judgment are:
\begin{mathpar}
   \Gamma \der_X \and %
   \Gamma \der_X A \and %
   \Gamma \der_X A = B  \and %
   \Gamma \der_X t : A  \and %
   \Gamma \der_X t = u : A  \and %
   \sigma \co \Delta \to_X \Gamma
\end{mathpar}
where the involved expressions are at stage $X$.
\begin{rem}
  In cubical type theory as described in~\cite{CCHM} we did not index
  judgments by objects $X$ but allowed extending context by interval
  variables instead.  Loosely speaking, a judgment $\Gamma
  \der_{\braces{i_1,\dots,i_n}} \JJ$ corresponds to $i_1 : \II, \dots,
  i_n : \II, \Gamma \der \JJ$ given the setting of~\cite{CCHM}.
\end{rem}
As mentioned above
we have the rule:
\begin{mathpar}
  \inferrule {\Gamma \der_X \JJ \\ f \co Y \to X} {\Gamma f \der_Y \JJ f}
\end{mathpar}

At each stage we have all the usual rules valid in a cwf with
$\Pi$-types, $\Sigma$-types, universes, and natural numbers.  We will
present some of the rules, but skip all congruence rules.
\begin{mathparpagebreakable}
  % basic cwf
  % contexts
  \inferrule{ }{\emptyctx\der_X} \and
  \inferrule{\Gamma\der_X \\ \Gamma\der_X A}{\Gamma.A\der_X} %
  \and%
  % types
  \inferrule %
  {\Gamma\der_X A \\ \sigma\co\Delta\to_X\Gamma}%
  {\Delta\der_X A\sigma}%
  \and
  %\\
  % terms
  \inferrule %
  {\Gamma\der_X t:A \\ \sigma \co \Delta\to_X\Gamma}%
  {\Delta\der_X t\sigma:A\sigma}%
  \and%
  \inferrule{\Gamma\der_X A}{\Gamma.A\der_X \qq:A\pp}%
  \and%
  \inferrule{\Gamma \der_X t : A \\ \Gamma \der_X A = B} {\Gamma
    \der_X t : B} %
  \and %
  %\\
  % substitutions
  \inferrule {\Gamma\der_X}{\id \co \Gamma \to_X \Gamma} %
  \and %
  \inferrule{\Gamma \der_X}{\emptysubst \co \Gamma \to \emptyctx}%
  \and %
  \inferrule{\Gamma\der_X A}{\pp \co \Gamma.A\to_X\Gamma}%
  \and%
  \inferrule %
  {\sigma \co \Delta\to_X \Gamma \\ \tau \co \Theta\to_X \Delta}%
  {\sigma\tau\co\Theta\to_X \Gamma}%
  \and%
  \inferrule{\sigma \co \Delta\to_X\Gamma \\ \Gamma\der_X A \\
    \Delta\der_X u:A\sigma}%
  {(\sigma,u)\co\Delta\to_X\Gamma.A}
  %
  % \\
  % \inferrule{\Gamma \der_X t = u : A \\ \Gamma \der_X A = B} {\Gamma
  %   \der_X t = u : B}
  \\
  % Pi types
  \inferrule{\Gamma.A\der_X B}{\Gamma\der_X\iPi(A,B)}%
  \and%
  \inferrule{\Gamma.A\der_X B \\ \Gamma.A\der_X b:B} %
  {\Gamma\der_X\iabs (b):\iPi(A,B)} %
  \and%
  \inferrule{\Gamma\der_X w:\iPi(A,B)\\\Gamma\der_X u:A} %
  {\Gamma\der_X \iapp(w,u):B[u]}
  % \\
  % % Sigma types
  % \inferrule{\Gamma.A\der_X B}{\Gamma\der_X\iSigma(A,B)}%
  % \and %
  % \inferrule{\Gamma.A\der_X B \\\Gamma\der_X u:A\\\Gamma\der_X v:B[u]} %
  % {\Gamma\der_X \ipair(u,v):\iSigma(A,B)} %
  % \and %
  % \inferrule{\Gamma\der_X w:\iSigma(A,B)}{\Gamma\der_X \ifst(w) : A}\and
  % \inferrule{\Gamma\der_X w:\iSigma(A,B)}{\Gamma\der_X \isnd(w) : B[\ifst(w)]}
\end{mathparpagebreakable}
where we write $[u]$ for $(\id,u)$ and $\sigma^+$ for $(\sigma \pp,
\qq)$.  The judgmental equalities (skipping suitable premises, types,
and contexts) are:
\begin{mathparpagebreakable}
  % equality rules
  %\\
  \id \, \sigma = \sigma \, \id = \sigma \and %
  (\sigma\tau)\delta = \sigma(\tau\delta) \and %
  \emptysubst \sigma = \emptysubst %
  \and
  (\sigma,u)\delta = (\sigma\delta,u\delta) \and \pp (\sigma,u) =
    \sigma \and \qq(\sigma,u) = u \and (\pp, \qq) = \id
  \and
  A \, \id = A \and (A \sigma) \delta = A (\sigma \delta) \and u \,
    \id = u \and (u \sigma) \delta = u (\sigma \delta)
  \\
  (\iPi(A,B))\sigma = \iPi(A\sigma,B \sigma^+) \and (\iabs
    (b)) \sigma = \iabs (b \sigma^+)
  \and
  \iapp(w,u)\delta= \iapp(w\delta,u\delta)\and \iapp(\iabs (b),u)=
  b[u]\and w = \iabs (\iapp(w \pp, \qq))
\end{mathparpagebreakable}
We skip the rules for $\Sigma$-types and natural numbers as they are
standard, but simply indexed with an object $X$ as we did for
$\Pi$-types.  The rules for universes are:
\begin{mathparpagebreakable}
  % universes
  \inferrule{\Gamma \der_X} {\Gamma \der_X \iU_n} %
  \and %
  \inferrule{\Gamma \der_X} {\Gamma \der_X \iU_n : \iU_{n+1}} %
  \and %
  \inferrule{\Gamma \der_X A : \iU_n} {\Gamma \der_X A : \iU_{n+1}} %
  \and %
  % \inferrule{\Gamma \der_X A = B : \iU_n} {\Gamma \der_X A = B :
  %   \iU_{n+1}} %
  % \and %
  \inferrule{\Gamma \der_X A : \iU_n} {\Gamma \der_X A } %
  % \and %
  % \inferrule{\Gamma \der_X A = B : \iU_n} {\Gamma \der_X A = B }
\end{mathparpagebreakable}
and we skip the rules for equality and closure under the type formers
$\iPi$,$\iSigma$, natural numbers, $\iPath$, and $\iGlue$.

To state the rules for dependent path-types we introduce the following
abbreviations.  We write $\Gamma.\II \der_X \bar A$ if $\bar A =
(A_{f,r})$ is a family indexed by $Y$, $f\co Y \to X$, and $r \in
\II(Y)$ such that
\[
  \Gamma f \der_Y A_{f,r} \text{ and } \Gamma f g\der_Z (A_{f,r})g =
  A_{fg,r}.
\]
Given $\Gamma.\II \der_X \bar A$ we write $\Gamma.\II\der_X \bar u :
\bar A$ whenever $\bar u = (u_{f,r})$ is a family indexed by $Y$,
$f\co Y \to X$, and $r \in \II(Y)$ such that
\[
  \Gamma f \der_Y u_{f,r} : A_{f,r} \text{ and } %
  \Gamma f g \der_Z (u_{f,r})g = u_{fg,rg} : A_{fg,rg}.
\]
The rules for the dependent path type are:
\begin{mathpar}
  \inferrule%
  {
    % \Gamma \der_Y A_{f,r} \text{ for each }f : Y \to X, r \in \II(Y) \\\\
    % \Gamma \der_Z (A_{f,r})g = A_{fg, rg} \text{ for each }f \co Y \to X, r \in \II(Y),
    % g \co Z \to Y \\\\
    \Gamma.\II \der_X \bar A\\
    \Gamma \der_X u: A_{\id_X,0}\\
    \Gamma \der_X u: A_{\id_X,1}
  }
  { \Gamma \der_X \iPath(\bar A,u,v) } %
  \and %
  \inferrule%
  {
    \Gamma.\II \der_X \bar A\\
    \Gamma .\II \der_X \bar u : \bar A \\
  }
  {
    \Gamma \der_X \iabs (\bar u) : \iPath(\bar A, u_{\id_X,0}, u_{\id_X,1})
  }
  \and %
  \inferrule%
  { \Gamma \der_X t : \iPath(\bar A,u,v) \\ r \in \II(X) }
  { \Gamma \der_X \ipapp(t,r) : A_{\id_X,r}}
  \\
  \ipapp(\iabs (\bar u), r) = u_{\id,r}  \and %
  t = \iabs (\ipapp(t f, r)_{f,r}) \and %
  \iPath (\bar A, u, v) \sigma = \iPath ((A_{f,r} \sigma f)_{f,r}, u \sigma, v
  \sigma) \and %
  (\iabs(\bar u)) \sigma = \iabs ((u_{f,r} \sigma f)_{f,r}) \and %
  (\ipapp (t,r)) \sigma = \ipapp (t \sigma, r)
\end{mathpar}
Note that in general these rules might have infinitely many premises.
We get the non-dependent path type for $\Gamma \der_X A$ by using the
family $A_{f,r} := Af$.

% Filling structure

Given $\Gamma.\II \der_X \bar A$ and $b \in \braces{0,1}$ we write
$\Gamma.\II \der_X^{\phi,b} \bar u : \bar A$ for $\bar u = (u_{f,r})$ a
family indexed over all $Y$, $f \co Y \to X$, and $r \in \II(Y)$ such
that either $f$ is in the sieve $[\phi]$ or $r = b$ and we have
\[
  \Gamma f \der_Y u_{f,r} : A_{f,r} \text{ and } \Gamma f g \der_Z
  (u_{f,r})g = u_{fg,rg} : A_{fg,rg}
\]
for all $g \co Z \to Y$.  The rule for the filling operation is given
by:
\begin{mathpar}
  \inferrule %
  { \Gamma.\II \der_X \bar A \\
    \phi \in \FF (X) \\
    b \in \braces{0,1} \\
    % \Gamma f \der_Y u_{f,r} : A_{f,r} \text{ for all } f \co Y \to X
    % \text{ and } r \in \II(Y) \text { such that } f \in \phi(Y) \text{
    %   whenever } r \neq b\\
    \Gamma.\II \der_X^{\phi,b} \bar u : \bar A\\
    r \in \II(X)
  }
  {
    \Gamma \der_Y \ifill(\bar A,\phi,b,\bar u,r) : A_{\id,r}
  } %
  % \and %
  % \ifill(A,\phi,b,u,r) = u_{\id,r} \text { whenever } \phi(X) = \top
  % \text{ or } r= b
\end{mathpar}
with judgmental equality
\begin{mathpar}
  \ifill(\bar A,\phi,b,\bar u,r) = u_{\id,r} \text { whenever } [\phi]
  \text{ is the maximal sieve or } r= b.
\end{mathpar}

For the glueing operation we only present the formation rule; the
other rules are similar as in~\cite{CCHM} but adapted to our setting.
We write $\Gamma \der_X^\phi \bar B$ if $\bar B$ is a family of $B_f$
for $f \co Y \to X$ in $[\phi]$ with $\Gamma f \der_Y B_f$ which is
compatible, \ie $\Gamma f g \der_Z B_{f}g = B_{fg}$.  In this case, we
write likewise $\Gamma \der_X \bar u : \bar B$ if $\bar u$ is a
compatible family of terms $\Gamma f \der_Y u_f : B_f$.
\begin{mathpar}
  \inferrule { %
    \Gamma \der_X A \\
    \phi \in \FF (X)\\
    \Gamma \der_X^\phi \bar B\\
    \Gamma \der_X^\phi \bar u : \iisEquiv (\bar B,A)\\
  }
  {\Gamma \der_X \iGlue (A,\phi,\bar B, \bar u)}
\end{mathpar}
and the judgmental equality $\iGlue (A,\phi,\bar B, \bar u) = B_\id$
in case $[\phi]$ is the maximal sieve, and an equation for
substitution.

\medskip

This formal system gives rise to a cubical cwf $\Syn$ as follows.
First, define judgmental equality for contexts and substitutions as
usual (we could also have those as primitive judgments).  Next, we
define presheaves $\Con$ and $\Subst$ on $\mathcal{C}$ by taking, say,
$\Con (X)$ equivalence classes $\eqcl \Gamma$ of $\Gamma$ with $\Gamma
\der_X$ modulo judgmental equality; restrictions are induced by the
(implicit) substitution: $\eqcl{\Gamma} f = \eqcl{\Gamma f}$.  Types
$\Type(X,\eqcl{\Gamma})$ are equivalence classes of $A$ with $\Gamma
\der_X A$ modulo judgmental equality, and elements are defined
similarly as equivalence classes.

For type formers in $\Syn$ let us look at path types: we have to give
an element of $\Type(\Gamma)$ in a context (w.r.t.\ the internal
language) $\Gamma : \Con, A : \II \to \Type(\Gamma), u :
\Elem(\Gamma,A\,0), v : \Elem(\Gamma, A\,1)$.  Unfolding the use of
internal language, given $\eqcl{\Gamma} \in \Con(X)$, a compatible
family $\eqcl{A_{f,r}} \in \Type(Y,\eqcl{\Gamma} f)$ (for $f \co Y \to
X$ and $r \in \II(Y)$) and elements $\eqcl{u} \in
\Elem(\eqcl{\Gamma},\eqcl{A_{\id,0}})$ and $\eqcl{v} \in
(\eqcl{\Gamma},\eqcl{A_{\id,1}})$, we have to give an element of
$\Type (X,\eqcl{\Gamma})$, which we do by the formation rule for
$\iPath$.

The remainder of the cubical cwf structure for $\Syn$ is defined in a
similar manner, in fact the rules are designed to reflect the laws of
cubical cwfs. We conjecture that we can follow a similar argument as in~\cite{Streicher91}
to show that $\Syn$ is the \emph{initial} cubical
cwf.  Given a cubical cwf $\MM$ over $\mathcal{C},\II,\FF$ we first
have to define \emph{partial} interpretations of the raw syntax and
then show that each derivable judgment has a defined interpretation in
$\MM$, and for equality judgments both sides of the equation have a
defined interpretation in $\MM$ and are equal.  In an intuitionistic framework,
this partial interpretation should be described as an inductively defined relation,
which is shown to be functional. The partial
interpretation $\brackss{ - }$ assigns meanings to raw judgments with
the following signature:
\begin{equation*}
  \begin{array}{lll}
    \brackss{\Gamma \der_X}
    & \in & \Con_\MM (X)
    \\
    \brackss{\sigma : \Delta \to \Gamma}
    & \in & \Subst_\MM (X,\brackss{\Delta \der_X},\brackss{\Gamma \der_X})
    \\
    \brackss{\Gamma \der_X A}
    & \in & \Type_\MM (X, \brackss{\Gamma \der_X})
    \\
    \brackss{\Gamma \der_X u : A}
    & \in & \Elem_\MM (X,  \brackss{\Gamma \der_X}, \brackss{\Gamma \der_X A})
  \end{array}
\end{equation*}
where among the conditions for the interpretation on the left-hand
side to be defined is that all references to the interpretation on the
right-hand side are defined.  This proceeds by structural induction on
the raw syntax and for $\brackss{\Gamma \der_X \JJ}$ to be defined we
assume all the ingredients needed are already defined.  \Eg for the
path type $\brackss{\Gamma \der_X \iPath (\bar A,u,v)}$ we in
particular have to assume that the assignment $f, r \mapsto
\brackss{\Gamma f \der_Y A_{f,r}}$ is defined and gives rise to a
suitable input of $\iPath_\MM$.

\section{Indexed inductive sets in presheaves over \texorpdfstring{$\mathcal{C}$}{C}}
\label{indexed-inductive-sets}

We work in the setting of \cref{developments} given by presheaves over $\mathcal{C}$.

Given a set $I$, a family $A$ over $I$, a family $B$ over $i : I$, an element $a : A\,i$, and a map
\[
s : \prod{i : I}\prod{a : A\,i} B\,i\,a \to I
,\]
the \emph{indexed inductive set} $\W_{I, A, B, s}$ is the initial algebra of the polynomial endofunctor~\cite{gambino-kock:polymomial} on the (internal) category of families over $I$ sending a family $X$ to the family
\[
\brackss{I, A, B, s}\,i = \sum{a : A\,i} \prod{b : B\,i\,a} X(s\,i\,a\,b)
.\]
Its constructive justification as an operation in the internal language of the presheaf topos using inductive constructions of the metatheory is folklore (in a classical setting, one can use transfinite colimits~\cite{kelly:transfinite-colimit}).%
\footnote{An indexed inductive set in presheaves unfolds externally to an indexed inductive-recursive definition where one defines the values at every level simultaneously with the restriction operations between levels.
In turn, this indexed inductive-recursive definition can be encoded as an indexed inductive-inductive definition, which in turn reduces to an indexed inductive definition (which one may further reduce to an inductive definition~\cite{gambino-hyland:inductive-families}).
Both steps use the idea of encoding functions via their graphs.
Alternatively, one can directly transform the inductive-recursive definition to an indexed inductive definition by first omitting about the naturality condition that mentions restriction, then define restriction recursively, and finally carve out the elements that recursively satisfy the naturality condition.}
If $I, A, B$ are small with respect to a universe $\U_i$ with $i \in \braces{0, 1, \ldots, \omega}$, then $\W_{I, A, B, s} : I \to \U_i$.

Let $I, A, B$ now be small with respect to $\U_\omega$.
Given $\Fill(I, A)$ and $\Fill(\sum{i : I}A\,i, \lam{(i, a)} B\,i\,a)$, we may use induction (\ie the universal property of $W_{I, A, B, s}$) to derive an element of $\Fill(I, \W_{I, A, B, s})$.
As in \cref{developments} for dependent products, this implies (using external reasoning) the internal statement $\CC(\W_{I, A, B, s}\,i)$ for $i : I$ given $\CC(A\,i)$ for all $i$ and $\CC(B\,i\,a)$ for all $i, a$.
We then call $W_{I, A, B, s}$ a \emph{fibrant uniformly indexed inductive set}.
The qualifier \emph{uniformly indexed} indicates that $A$ is a fibrant family over $I$ rather than a fibrant set with a ``target'' map to $I$ that indicates the target sort of the constructor $\sup$.

Given $A : \U_\omega$ with $\Fill(1, A)$, we may use the technique of Andew Swan~\cite{swan:identity-type,OP} to construct a (level preserving) identity set $\Id_A\,a_0\,a_1$ for $a_0, a_1 : A$ (different from the equality set $a_0 = a_1$) with $\Fill(A \times A, \lam{(a_0, a_1)} \Id_A\,a_0\,a_1))$ and constructor $\refl_a : \Id_A\,a\,a$ for $a : A$ that has the usual elimination with respect to families $P : \prod{a_0\,a_1 : A} \Id_A\,a_0\,a_1 \to \U_\omega$ that satisfy
$\Fill(\sum{a_0\,a_1 : A} \Id_A\,a_0\,a_1, P)$.
Using external reasoning as before, one has $\CC(\Id_A\,a_0\,a_1)$ given $\CC(A)$, justifying calling $\Id_A\,a_0\,a_1$ a \emph{fibrant identity set}; using~\eqref{U-fib-has-fill} one has elimination with respect to families $P$ of the previous signature with $\CC(C\,a_0\,a_1\,p)$ for all $a_0, a_1, p$.

Using a folklore technique, we may use fibrant identity sets to derive \emph{fibrant indexed inductive sets} from fibrant uniformly indexed inductive sets, by which we mean the following.
Given $(I, \fib_I) : \Ufib_\omega$, $(A, \fib_A) : \Ufib_\omega$, $\angles{B, \fib_B} : A \to \Ufib_\omega$ with maps $t : A \to I$ and $s : \prod{a : A} B\,a \to I$, we have $\angles{W_{I, A, B, s, t}, \fib_W} : I \to \Ufib_\omega$ (we omit the subscripts to $W$ for readability), $W$ living in $\U_i$ if $I, A, B$ do,  with
\[
\sup : \prod{a : A}\prod{f : \prod{b : B\,a} W\,(s\,a\,b)} W\,(t\,a)
.\]
Given $\angles{P, \fib_P} : \prod{i : I} W \to \Ufib_\omega$ with
\[
h : \prod{a : A}\prod{f : \prod{b : B\,a} W\,(s\,a\,b)} (\prod{b : B\,a} P\,(s\,a\,b)\,(f\,b)) \to P\,(t\,a)\,(\sup\,a\,f)
,\]
we have $v : \prod{i : I}\prod{w : W\,i} P\,i\,w$ such that
\[
v\,(t\,a)\,(\sup\,a\,f) = h\,a\,f\,(\lam{b} v\,(s\,a\,b)\,(f\,b)
.\]

Fibrant indexed inductive sets are used for the interpretation in the sconing model of natural numbers in \cref{sconing}, higher inductive types in \cref{higher-inductive-types}, and identity types in \cref{identity-types}.
In practise, we will usually not bother to bring the fibrant indexed inductive set needed into the above form and instead work explicitly with the more usual specification in terms of a list of constructors, each taking a certain number non-recursive and recursive arguments.%
\footnote{
Note that the latter is really an instance of the former since our dependent sums, dependent products, and finite coproducts are extensional (satisfy universal properties).
Conversely, the former is an instance of the latter with a single constructor taking a non-recursive and a recursive argument.
}

As an example, we construct the fibrant indexed inductive set $\iN'$ needed in \cref{sconing}.
There, we have a fibrant set $\verts{\iN} : \U_0$ (satisfying $\CC(\verts{\iN})$) with an element $\izero : \verts{\iN}$ and an endofunction $\isucc : \verts{\iN} \to \verts{\iN}$.
We wish to define the fibrant indexed inductive set $\iN' : \verts{\iN} \to \U_0$ with constructors $\zero' : \iN'\,\izero$ and $\succ' : \prod{n : \verts{\iN}\,\rho} \iN'\,n \to \iN'\,(\isucc\,n)$.
We let $\iN'$ be the uniformly indexed inductive set over $m : \verts{\iN}$ with constructors
\begin{align*}
\zero'' &: \Id_{\verts{\iN}}\,m\,\izero \to \iN'\,m
,\\
\succ'' &: \prod{n : \verts{\iN}} \Id_{\verts{\iN}}\,m\,(\isucc\,n) \to \iN'\,n \to \iN'\,m
.\end{align*}
and define $\zero' = \zero''\,\refl_\izero$ and $\succ'\,n\,n' = \succ''\,n\,\refl_{\isucc(n)}\,n'$.
Fibrancy of $\Id$ ensures fibrancy of $\iN'$ (\ie $\CC(\iN'\,n)$ for $n : \verts{\iN}$).
For elimination, we are given a fibrant family $P\,n\,n'$ for $n : \verts{\iN}$ and $n' : \iN'\,n$ with $z' : P\,\izero\,\izero'$ and $s'\,n\,n'\,x : P\,(\isucc\,n)\,(\isucc'\,n\,n')$ for all $n, n'$ and $x : P\,n\,n'$.
We have to define $h'\,n\,n' : P\,n\,n'$ for all $n, n'$ such that $h'\,\izero\,\izero' = z'$ and $h'\,(\isucc\,n)\,(\isucc'\,n\,n') = s'\,n\,n'\,(h'\,n\,n')$.
We define $h'$ by induction on the uniformly indexed inductive set $\iN'$ and fibrant identity sets (using fibrancy of $P$) via defining equations
\begin{align*}
h'\,\zero''\,\refl_\izero &= z
,\\
h'\,(\isucc n)\,(\succ''\,n\,\refl_{\isucc\,n}\,n') &= s'\,n\,n'\,(h'\,n\,n')
.\end{align*}

\section{Variations}
\label{appendix-variations}

\subsection{Univalence as an axiom}
\label{univalence-as-axiom}

Our treatment extends to the case where the glue types in a cubical cwf as in \cref{cubical-cwf} are replaced by an operation $\Elem(\Gamma, \iUnivalence_n)$ for $\Gamma : \Con$ and $n \geq 0$, with $\iUnivalence_n$ defined in \cref{univalent-type-theory-prop-J}.

To define this operation in the sconing model of \cref{sconing}, one first shows analogously to \cref{efun-pres-equiv,equiv-in-sconing} that $\verts{-}$ preserves contractible types and that $(A, A') : \Type^*(\Gamma, \Gamma')$ is contractible exactly if $A$ is contractible and $A'\,\rho\,\rho'\,a$ for $\rho : \verts{\Gamma}$ and $\rho' : \Gamma'\,\rho$ where $a : \verts{A}$ is the induced center of contraction.
We have analogous statements for types of homotopy level $n \geq 0$i in $\MM$, in which case we instead have to quantify over all $a : \verts{A}$.

Given $(A, A') : \Type_n^*(\Gamma, \Gamma')$, we have show that the type
\[
(S, S') = \iSigma^*(\iU_n, \Equiv^*(\qq, A))
\]
over $(\Gamma, \Gamma')$ is contractible in $\MM^*$.
Without loss of generality, we may assume the center of contraction of univalence in $\MM$ is given by the identity equivalence.
Using the observations of the preceding paragraph, it suffices to show that
\[
V' = S'\,\rho\,\rho'\,(\ipair(A \rho, (\iabs(\qq), w)))
\]
is contractible for $\rho : \verts{\Gamma}$ and $\rho' : \Gamma'\,\rho$ where $w$ denotes the canonical witness that the identity map $\iabs(\qq)$ on $A \rho$ is an equivalence in $\MM$.
Inhabitation is evident, and so it remains to show propositionality.
By the case of the preceding paragraph for propositions, the second component of $V'$ is a proposition, and thus we can ignore it for the current goal, which then becomes
\[
\isProp \parens[\big]{\sum{T' : \verts{A} \to \U_n} \prod{a : \verts{A}} \Equiv(T'\,a, A'\,\rho\,\rho'\,a)}
\]
and follows from univalence in the standard model, justified by glueing.

\section{Simplicial set model}
\label{simplicial-sets}

Choosing for $\mathcal{C}$ the simplex category $\Delta$, for $\II$ the usual interval $\Delta^1$ in simplicial sets, and for $\FF$ a small copy $\Omega_{0,\dec}$ of the sublattice of $\Omega_0$ of decidable sieves, we obtain a notion of cubical cwf with a simplicial notion of shape.

Assume now the law of excluded middle.
The above choice of $\mathcal{C}, \II, \FF$ satisfies all of the assumptions of \cref{developments} but one: the existence of a right adjoint to exponentiation with $\II$.
However, except for \cref{strict-canonicity}, the only place our development makes use of this assumption is in establishing~\eqref{C-vs-fill-external}.
We will instead give a different definition of $\CC$ that still satisfies~\eqref{C-vs-fill-external}.
Then the rest of our development, except for \cref{strict-canonicity}, still applies to simplicial sets.

A \emph{Kan fibration structure} on a family $Y \co X \to \U_\omega$ in simplicial sets consists of a choice of diagonal fillers in all commuting squares of the form
\begin{equation*} \label{kan-fibration-structure}
\begin{gathered}
\xymatrix{
  \Lambda_k^m
  \ar[r]
  \ar[d]
&
  \sum{x:X} Y\,x
  \ar[d]
\\
  \Delta^m
  \ar[r]
  \ar@{.>}[ur]
&
  X
}
\end{gathered}
\end{equation*}
with left map a \emph{horn inclusion} and right map the evident projection.
Note that the codomains of horn inclusions are representable.
It follows that the presheaf of Kan fibration structures indexed over the slice of simplicial sets over $\U_\omega$ is representable.
Given $[n] \in \Delta$ and $A \in (\U_\omega)_n$ (\ie an $\omega$-small presheaf on $\Delta/[n]$), we define $\CC([n], A)$ as the set of Kan fibration structures on $A \co \Delta^n \to \U_\omega$.
This defines a level preserving map $\CC \co \U_\omega \to \U_\omega$.
Then the representing object of the above presheaf is given by the first projection $\Ufib_\omega \to \U_\omega$ where $\Ufib_\omega = \sum{X:\U_\omega} \CC(X)$ is defined as before.

Let us now verify~\eqref{C-vs-fill-external}.
Given a simplicial set $X$ with $Y \co X \to \U_\omega$, a global element of $\Fill(X, Y)$ corresponds to a \emph{uniform Kan fibration} structure on $\sum{x:X}(Y\,x) \to X$ in the sense of~\cite{gambino-sattler:frobenius}.
A uniform Kan fibration structure induces a Kan fibration structure naturally in $X$, giving the forward direction of~\eqref{C-vs-fill-external}.
For the reverse direction, it suffices to give a uniform Kan fibration structure in the generic case, \ie a global element of $\Fill(\Ufib_\omega, \lam{(A, c)} A)$.
This is~\cite[Theorem~8.9, part~(ii)]{gambino-sattler:frobenius} together with the fact proved in~\cite[Chapter~IV]{gabriel-zisman:calculus-of-fractions} that Kan fibrations lift against pushout products of interval endpoint inclusions with (levelwise decidable) monomorphisms.%
\footnote{This is the only place where excluded middle is used, to produce a cellular decomposition in terms of simplex boundary inclusions of such a monomorphism.}

Having verified~\eqref{C-vs-fill-external}, the rest of our development applies just as well to the case of simplicial sets.
In particular, we obtain in the standard model $\Std$ of \cref{standard-model} a version of the simplicial set model~\cite{simplicial-set-model} of univalent type theory (using \cref{identity-types} for identity types).%
\footnote{Instead of Kan fibration structures, we can also work with the property of being a Kan fibration.
Then $\CC$ is valued in propositions and we would obtain in $\Std$ a version of the simplicial set model in which being a type is truly just a property.
However, choice would be needed to obtain~\eqref{C-vs-fill-external}.}
As per \cref{higher-inductive-types}, we furthermore obtain higher inductive types in the simplicial set model in a way that avoids (as suggested by Andrew Swan~\cite{swan:higher-inductive-types}) the pitfall of fibrant replacement failing to preserve size encountered in~\cite{lumsdaine-shulman:hits}.

Seeing simplicial sets as a full subtopos of distributive lattice cubical sets as observed in~\cite{kapulkin:cubical}, there is a functor from cubical cwfs with $(\mathcal{C}, \II, \FF) = (\Delta, \Delta^1, \Omega_{0,\dec})$ to cubical cwfs where $\mathcal{C}$ is the Lawvere theory of distributive lattices, $\II$ is represented by the generic object, and $\FF$ is the (small) sublattice of $\Omega_0$ generated by distributive lattice equations.
The cubical cwfs in the image of this functor satisfy a \emph{sheaf condition}, which can be represented syntactically as an operation allowing one to \eg uniquely glue together to a type $\Gamma \vdash_{\braces{i, j}} A$ coherent families of types $\Gamma f \vdash_X A_f$ for $f$ a map to $X$ from the free distributive lattice on symbols $\braces{i, j}$ such that $f\,i \leq f\,j$ or $f\,j \leq f\,i$ (compare also the tope logic of~\cite{riehl-shulman:synthetic-higher-cats}).

Applying this functor to the simplicial set model $\Std$ discussed above, we obtain an interpretation of distributive lattice cubical type theory (with $\II$ and $\FF$ as above) in the sense of the current article (crucially, without computation rules for filling at type formers) in simplicial sets.
Thus, this cubical type theory is homotopically sound: can only derive statements which hold for standard homotopy types.

\vspace{2cm}\phantom{X}

% Fix for lmcs footer causing footnote to overlap with page text.
% \newpage

\end{document}